\documentclass[preprint]{elsarticle}

\usepackage{amsmath,amssymb,bm}
\usepackage{graphicx}
\usepackage{caption,subcaption}
\usepackage{geometry}
\usepackage{booktabs,siunitx,array}
\usepackage{algorithm}
\usepackage{algpseudocode}
\usepackage{hyperref}
\usepackage{cleveref}
\usepackage{enumerate}
\usepackage{float}
\usepackage{placeins}
\usepackage{color}
\usepackage{lineno}
\usepackage{makecell}
\usepackage{hyperref}       % hyperlinks
\usepackage{url}            % simple URL typesetting
\usepackage{amsfonts}       % blackboard math symbols
\usepackage{nicefrac}       % compact symbols for 1/2, etc.
\usepackage{microtype}      % microtypography
\usepackage{natbib}

\usepackage{doi}
\usepackage{dcolumn}% Align table columns on decimal point
\usepackage{bm}% bold math
\usepackage{caption}
\usepackage{subcaption}
\usepackage{float}
\usepackage{array}
\newcolumntype{P}[1]{>{\centering\arraybackslash}p{#1}}

\usepackage{xcolor}

\date{}
\begin{document}

\begin{frontmatter}

\title{Differentiable Inverse Modeling with Physics-Constrained Latent Diffusion for Heterogeneous Subsurface Parameter Fields}

\author[UMN]{Zihan Lin}
\author[UMN]{QiZhi He\corref{mycorrespondingauthor}}
\cortext[mycorrespondingauthor]{Corresponding author}
\ead{qzhe@umn.edu}
\address[UMN]{Department of Civil, Environmental, and Geo- Engineering, University of Minnesota, 500 Pillsbury Drive S.E., Minneapolis, MN 55455}
 
\begin{abstract}
\noindent We present a latent diffusion–based differentiable inversion method (LD-DIM) for PDE-constrained inverse problems involving high-dimensional spatially distributed coefficients. 
LD-DIM couples a pretrained latent diffusion prior with an end-to-end differentiable numerical solver to reconstruct unknown heterogeneous parameter fields in a low-dimensional nonlinear manifold, improving numerical conditioning and enabling stable gradient-based optimization under sparse observations.
The proposed framework integrates a latent diffusion model (LDM), trained in a compact latent space, with a differentiable finite-volume discretization of the forward PDE. Sensitivities are propagated through the discretization using adjoint-based gradients combined with reverse-mode automatic differentiation. 
Inversion is performed directly in latent space, which implicitly suppresses ill-conditioned degrees of freedom while preserving dominant structural modes, including sharp material interfaces.
The effectiveness of LD-DIM is demonstrated 
using a representative inverse problem for flow in porous media,
where heterogeneous conductivity fields are reconstructed from spatially sparse hydraulic head measurements. 
Numerical experiments assess convergence behavior and reconstruction quality for both Gaussian random fields and bimaterial coefficient distributions.
The results show that LD-DIM achieves consistently improved numerical stability and reconstruction accuracy of both parameter fields and corresponding PDE solutions compared with physics-informed neural networks (PINNs) and physics-embedded variational autoencoder (VAE) baselines, while maintaining sharp discontinuities and reducing sensitivity to initialization. 
% This indicates that LD-DIM provides a scalable and fully differentiable computational paradigm for integrating deep learning generative priors with classical PDE-constrained optimization in large-scale inverse problems.
\end{abstract}

% --------------------------------------------------------------------
\end{frontmatter}
\section{Introduction}\label{sec:intro}
Inverse problems constrained by partial differential equations (PDEs) arise ubiquitously in physics and engineering, where unknown spatially distributed parameters need to be inferred from sparse and indirect observations. 
Such problems are typically ill-posed and high-dimensional, exhibiting severe non-uniqueness and strong sensitivity to data availability, parameterization, and numerical conditioning,
particularly when the governing physics involves heterogeneous media and sharp spatial contrasts.
% When the governing physics involves heterogeneous media and sharp spatial contrasts, these challenges are further amplified, rendering stable and efficient inversion computationally demanding.

A canonical example of such challenges arises in subsurface flow modeling, where hydraulic properties such as conductivity exert a dominant influence on flow and transport processes in groundwater systems~\cite{Brunner2012}. 
These properties cannot be measured exhaustively due to the prohibitive cost and spatial sparsity of borehole investigations, and must instead be inferred indirectly from hydraulic or tracer observations \citep{Freeze1979,hill2007effective,Carrera2005}. 
As a result, inverse modeling has become a central computational tool for subsurface characterization \citep{Yeh2000,ramm2005inverse,linde2015geological}. 
However, subsurface inverse problems are notoriously ill-posed \citep{zhou2014inverse,yeh1986review}, owing to the combination of high-dimensional parameter spaces \citep{mclaughlin1996reassessment}, limited and unevenly informative observations \citep{Neuman1973Calibration}, and spatially variable parameter sensitivities. 
% For example, hydraulic head measurements near prescribed head boundaries are often dominated by boundary conditions rather than local conductivities, significantly reducing their diagnostic value.
These difficulties are further exacerbated by strong geological heterogeneity, which introduces sharp spatial contrasts in subsurface properties and amplifies numerical ill-conditioning~\cite{golmohammadi2015group}.

A widely adopted strategy for mitigating ill-posedness is to reduce the effective dimensionality of the parameter space \citep{mclaughlin1996reassessment,cui2014likelihood}. 
Numerous low-dimensional parameterization methods have been proposed, including transform-based approaches such as the discrete cosine transform (DCT) \citep{jafarpour2008history}, statistical techniques such as principal component analysis (PCA) and its variants (K-PCA, O-PCA) \citep{sarma2008kernel,vo2015data,lee2014large}, and stochastic representations such as the Karhunen–Loève expansion (KLE) \cite{tartakovsky2021physics,yeung2024gaussian}.
While effective in certain settings, these methods are typically built on linear or Gaussian assumptions and rely primarily on second-order statistics. 
As a consequence, they struggle to represent nonlinear, non-Gaussian geological features, often producing overly smooth parameter fields that fail to preserve connectivity and sharp interfaces essential for realistic flow and transport behavior \citep{chan2019parametric,linde2015geological}.

Recent advances in deep generative modeling have opened new avenues for nonlinear dimensionality reduction in inverse problems. 
Autoencoders and variational autoencoders (VAEs) learn compact latent representations that capture complex spatial patterns beyond the reach of linear techniques \citep{hinton2006reducing,masci2011stacked}. 
In subsurface applications, Laloy et al.~\citep{laloy2017inversion} demonstrated that VAE-based parameterizations can achieve orders-of-magnitude dimensionality reduction while substantially accelerating Markov chain Monte Carlo (MCMC) inversion. 
Subsequent studies have integrated such latent representations with differentiable forward solvers to enable gradient-based inversion \citep{wu2023physics}. 
Nevertheless, VAE-based approaches often yield overly smooth reconstructions and exhibit limited ability to preserve sharp geological contrasts when observational data are sparse.
To preserve geological realism and multiple-point statistics, Laloy et al.~\citep{laloy2018training} proposed spatial generative adversarial networks (GANs) to learn low-dimensional representations. However, GAN-based approaches usually suffer from training instability and strong sensitivity to hyperparameter choices.

%% 2. DDPMs
More recently, generative diffusion models have emerged as a powerful and stable  paradigm for high-fidelity image synthesis~\citep{ho2020denoising,rombach2022high}. 
By progressively transforming noise into structured samples, diffusion models have demonstrated superior expressiveness and robustness compared to VAEs and GANs, and have shown strong potential as priors for inverse problems \citep{rombach2022high,wang2023prior}. 
%% 3. LDM
Latent diffusion models (LDMs) further enhance computational efficiency by performing diffusion in a compact latent space while retaining strong generative capacity~\cite{rombach2022high}. 
In the subsurface context, recent studies have shown that LDMs can learn geologically informed latent manifolds that encode multiscale heterogeneity and complex structural patterns \citep{wang2023prior,zhan2025toward,lee2025latent,di2024latent,di20253d}. 
For example, Zhan et al.~\citep{zhan2025toward} demonstrated that LDMs can learn cross-scale representations of subsurface heterogeneity and establish consistent mappings among geological structures, flow responses, and concentration fields within a unified framework. 
Related work has also explored conditional LDMs for subsurface modeling, including the generation of reservoir facies distributions and associated responses constrained by prescribed geological patterns \citep{lee2025latent}. 
These studies suggest that LDMs are capable of learning compact, geologically informed latent manifolds and well suited for inverse modeling.

Building on this direction, Di Federico and Durlofsky~\citep{di2024latent} proposed a physics-guided LDM for low-dimensional parameterization of facies-based geomodels and demonstrated successful history matching via ensemble smoother with multiple data assimilation (ESMDA) in two-dimensional channelized systems. 
However, most existing deep generative model–based inversion approaches \citep{di20253d,di2024latent,chen2023fracture,dasgupta2025conditional} usually rely on ensemble-based calibration, such as ESMDA \citep{emerick2012ensemble,di2024latent}, or Bayesian posterior sampling, which require numerous computationally prohibitive forward simulations. 
This computational burden limits scalability and motivates inversion frameworks that directly exploit differentiability and latent-space optimization for large-scale subsurface flow problems.

%%%%%%%%%%%%%%%%%%%%%%%%% Differentiable Programming
In parallel, advances in differentiable numerical simulation \citep{innes2019differentiable} have enabled discretized PDE solvers to be embedded directly within optimization and inverse modeling pipelines via automatic differentiation \citep{baydin2018automatic}. 
This paradigm has received increasing attention in computational geosciences \citep{shen2023differentiable,gelbrecht2023differentiable,wu2023physics}, solid mechanics \citep{xue2023jax,du2024neural,du2025jax}, and fluid mechanics \citep{bezgin2023jax}. 
By propagating gradients through the numerical discretization, differentiable solvers provide a flexible alternative to classical adjoint implementations and enable end-to-end gradient-based inversion.
It is worth noting that automatic differentiation of PDEs also forms a core component in the development of physics-informed neural networks (PINNs) \citep{raissi2019physics,karniadakis2021physics}. 
However, the key distinction between these two paradigms lies in how the governing equations are incorporated into the learning process: PINNs embed the continuous PDE formulation directly into the loss function, whereas differentiable numerical solvers differentiate through a fully discretized forward model within the optimization loop.

Motivated by these developments, this work proposes the \emph{latent Diffusion-based differentiable inversion method} (LD-DIM), a computational framework that couples a pretrained latent diffusion model with a fully differentiable numerical solver to reconstruct PDE-constrained inverse problems in a low-dimensional latent space. 
In LD-DIM, the LDM provides a nonlinear mapping from low-dimensional latent variables to spatially distributed parameter fields, acting as an implicit prior that suppresses ill-conditioned modes while preserving complex geological structures. 
The forward physics is solved using a differentiable finite-volume discretization of Darcy flow, and sensitivities with respect to latent variables are obtained via end-to-end reverse-mode automatic differentiation combined with adjoint-based gradient evaluation~\cite{dwight2006effect,plessix2006review,xiao2021deep}. 
This formulation enables direct gradient-based optimization in latent space, substantially improving numerical stability and computational efficiency relative to inversion in the original high-dimensional parameter space.

The inverse problem for flow in heterogeneous porous media (Darcy flow) is employed as a representative testbed to examine the computational behavior of the proposed framework, including convergence characteristics, sensitivity to observation sparsity, and the ability to recover both smooth Gaussian fields and sharp-interface bimaterial structures.
Numerical experiments demonstrate that LD-DIM outperforms existing approaches, including PINNs~\citep{raissi2019physics} applied to subsurface flow~\citep{he2020physics} and physics-embedded VAEs~\citep{wu2023physics}, achieving up to a fivefold reduction in reconstruction error while preserving geological discontinuities. 
To support efficient forward and backward computations, the finite-volume solver is implemented in \texttt{JAX}~\citep{frostig2019compiling,jax2018github}, enabling efficient reverse-mode automatic differentiation for gradient-based optimization.

The remainder of this paper is organized as follows. 
\Cref{sec:methodology} presents the proposed methodology, including latent diffusion model training (\Cref{sec:ldm}) and the LD-DIM framework (\Cref{sec:physics_inverse}). 
\Cref{sec:numerical-validation} provides comprehensive numerical validation through two case studies: Gaussian conductivity fields with varying correlation lengths (\Cref{subsec:gaussian-fields}) and complex bimaterial distributions with sharp interfaces (\Cref{subsec:bimaterial-fields}), along with performance comparisons against established methods. 
Finally, \Cref{sec:conclusion} concludes the paper with a discussion of key findings and future research directions.

\section{Methodology}\label{sec:methodology}
This section introduces the latent diffusion differentiable inversion model (LD-DIM), a physics-constrained inversion framework that integrates a pretrained latent diffusion model with a differentiable discretization-based solver for groundwater flow. The core idea is to encode the high-dimensional conductivity field in a compact latent space learned by the LDM, enabling efficient inversion over a significantly reduced set of unknowns using adjoint-based gradients while strictly enforcing the governing flow physics.

% --------------------------------------------------------------------

\subsection{Forward problem: flow in porous media}\label{sec:forward}
We consider steady-state groundwater flow in a confined, isotropic aquifer as a representative example. 
For single-phase flow in heterogeneous porous media, the governing equations consist of Darcy’s law:
 \begin{equation}\label{eq:darcy}
     \bm{u} = -K \nabla h, \quad \text{in } \Omega,
 \end{equation}
 where $\bm{u}$,  $K$, and $h$ denote the Darcy velocity,  hydraulic conductivity, and  pressure head, respectively,  
 and the continuity equation:
 \begin{equation}\label{eq:continuum}
    \nabla \cdot \bm{u} = 0, \quad \text{in } \Omega.
 \end{equation}
which applies in the absence of volumetric sources or sinks.
The spatial heterogeneity of the conductivity field is represented through a continuum description $K(\bm{x})$. 
The physical units are: $\bm{u}~[\mathrm{m/s}]$, $K~[\mathrm{m/s}]$, and $h~[\mathrm{m}]$.

Combining Eqs.~\eqref{eq:darcy} and \eqref{eq:continuum} yields the standard elliptic equation governing steady groundwater flow,
\begin{equation}\label{eq:darcy_governing}
   \nabla \cdot \bigl(K(\bm{x}) \nabla h(\bm{x}) \bigr) = 0 \quad \text{in } \Omega.
\end{equation}
The forward problem can be expressed compactly as the mapping $\mathcal{M}: K \mapsto  h$, where $\mathcal{M}$ denotes the solution operator that computes the hydraulic head field corresponding to a given conductivity distribution.

In this study, we develop a differentiable finite-volume solver to compute the steady-state solution of Eq.~\eqref{eq:darcy_governing}, as detailed in Section~\ref{sec:fvm}.

% --------------------------------------------------------------------

\subsubsection*{Challenges in conductivity estimation}
As discussed in the Introduction, estimating a heterogeneous conductivity field from sparse observations leads to a severely ill-posed inverse problem. The high dimensionality of the unknown field, together with the limited and spatially biased nature of available data, results in strong non-uniqueness and uncertainty.
To alleviate this ill-posedness, in this study we employ the generative latent diffusion model (LDM) that learns a compact latent representation of conductivity fields. 
This low-dimensional latent space reduces the effective number of unknowns while embedding prior geological structure, yielding a more tractable and physically consistent inversion framework.

% --------------------------------------------------------------------
\subsection{Latent diffusion model (LDM) as a pretrained surrogate}\label{sec:ldm}

LDMs offer a powerful parameterization framework by encoding high-dimensional hydrological fields into low-dimensional latent representations~\cite{di2024latent,zhan2025toward,zhang2024non}. 
As illustrated in Fig. \ref{fig:ldm_training}, an LDM consists of two major components: a variational autoencoder (VAE) that compresses input fields into latent variables, and a U-Net denoiser~\cite{ronneberger2015u} that performs the reverse diffusion process within this latent space.
In the remainder of this subsection, we review the basic components of LDMs for completeness. 
We first introduce classical diffusion models, followed by the VAE architecture used for latent compression. 
We then describe the training procedure and diffusion process of the full LDM within the latent space.
For additional details, we refer readers to~\cite{ho2020denoising,rombach2022high}.

\subsubsection{Diffusion model}\label{sec:diff_model}

Diffusion models are a class of generative models that transform data into noise through a forward diffusion process and recover it through a learned reverse denoising process. 
In the standard formulation, denoising diffusion probabilistic models (DDPMs)~\citep{ho2020denoising} add Gaussian noise to data via a Markov chain governed by a noise schedule $\{\beta_t\}_{t=1}^T$
and train a neural network to predict the injected noise at each timestep.
Denoising diffusion implicit models (DDIMs)~\citep{song2020denoising} extend DDPMs by introducing a deterministic reverse process that shares the same forward procedure and training objective. 
This deterministic formulation in DDIMs enables sampling with reduced reverse steps, which results in significantly faster inference, without requiring retraining, and is therefore adopted in this study.

\textit{I. Forward process}. The forward (diffusion) process follows a Markov chain that gradually adds Gaussian noise to the data according to a noise schedule $\{\beta_t\}_{t=1}^T$, where $T$ denotes the total number of diffusion steps. 
Define $\alpha_t:=1-\beta_t$ and $\bar{\alpha}_t:=\prod_{s=1}^t\alpha_s$. Each transition step is a Gaussian distribution:
\begin{equation}\label{eq:forward_transition}
q(z_t\mid z_{t-1})=\mathcal{N}\!\bigl(z_t;\sqrt{\alpha_t}\,z_{t-1},\,\beta_t\,\bm{I}\bigr).
\end{equation}
Accordingly, the marginal distribution at any timestep admits a closed-form expression:
%the distribution at any time step admits a closed-form expression:
\begin{equation}\label{eq:forward_marginal}
q(z_t\mid z_0)=\mathcal{N}\!\bigl(z_t;\sqrt{\bar{\alpha}_t}\,z_0,\,(1-\bar{\alpha}_t)\,\bm{I}\bigr),
\end{equation}
which can equivalently be sampled using the reparameterization
\begin{equation}\label{eq:forward_reparam}
z_t=\sqrt{\bar{\alpha}_t}\,z_0+\sqrt{1-\bar{\alpha}_t}\,\varepsilon,\qquad \varepsilon\sim\mathcal{N}(\mathbf{0},\bm{I}).
\end{equation}

% \textbf{DDIM reverse process.} 
\textit{II. Reverse process.}
The reverse process removes noise by predicting the noise component at each timestep. 
The learnable denoiser $\varepsilon_\theta(z_t,t)$, which is typically implemented with a U-Net architecture and parameterized by $\theta$ (see Fig. \ref{fig:ldm_training}b), is used together with the noise schedule $\{\alpha_t\}_{t=1}^T$.
%with $\bar\alpha_t=\prod_{s=1}^t\alpha_s$
The DDIM deterministic update step is given by 
\begin{equation}\label{eq:ddim_update}
z_{t-1}=
\sqrt{\alpha_{t-1}}\,
\frac{z_t-\sqrt{1-\bar\alpha_t}\,\varepsilon_\theta(z_t,t)}{\sqrt{\bar\alpha_t}}
\;+\;
\sqrt{1-\alpha_{t-1}}\;\varepsilon_\theta(z_t,t),
\end{equation}
for $t=T,\ldots,1$.
This deterministic denoising process constitutes the core generative mechanism of the DDIM-based diffusion model.

\textit{III Loss function.}
The neural network model $\varepsilon_\theta$ in Eq. \eqref{eq:ddim_update} is trained using the simplified noise-matching objective:

\begin{equation}\label{eq:diffusion_loss}
\mathcal{L}_{\text{diff}}(\theta) = \mathbb{E}_{t, z_0, \varepsilon}\left[\left\|\varepsilon - \varepsilon_\theta\left(\sqrt{\bar{\alpha}_t} z_0 + \sqrt{1 - \bar{\alpha}_t} \varepsilon, t\right)\right\|_2^2\right].
\end{equation}
where $\varepsilon\sim\mathcal{N}(\mathbf{0},\bm{I})$ is a standard Gaussian noise sample, which is also the ground truth noise to be predicted by the network.
This objective trains the U-Net denoiser to predict the noise injected at each forward diffusion timestep. 
The same trained model supports both DDPM stochastic sampling and DDIM deterministic generation.
In DDIM sampling, one initializes $z_T \sim \mathcal{N}(\mathbf{0}, \bm{I})$ and iteratively applies the update in Eq.~\eqref{eq:ddim_update} for $t = T, \ldots, 1$ to obtain $z_0$.

\subsubsection{Variational autoencoder}
To enable diffusion processes to operate in a compact domain, 
a VAE is employed  to learn a latent representation of conductivity fields.
As shown in Fig.~\ref{fig:ldm_training}a, the VAE consists of an \textit{encoder} $E_{\phi}$ and a \textit{decoder} $D_{\psi}$.

Given a conductivity field $K$, the encoder produces a pair
$(\bm{z}_{\mu},\bm{z}_{\sigma})=E_{\phi}(K)$,
which parameterize an approximate Gaussian posterior over latent variables conditioned on the data:
\begin{equation}\label{eq:vae_posterior}
  q_{\phi}(\bm{z}_0\mid K)=\mathcal{N}\!\bigl(\bm{z}_{\mu},\,\operatorname{diag}(\bm{z}_{\sigma}^{2})\bigr),
\end{equation}
where $\bm{z}_{\mu},\bm{z}_{\sigma}\in\mathbb{R}^{d}$ are the predicted mean and standard deviation vectors, and $d$ is the latent dimension. 
To enable backpropagation through the stochastic sampling process, a latent variable is drawn using the reparameterization trick:
\[\bm{z}_0=\bm{z}_{\mu}+\bm{z}_{\sigma}\odot\bm{\varepsilon},\qquad
  \bm{\varepsilon}\sim\mathcal{N}(\bm{0},\bm{I}_{d}),
\]
where $\bm{z}_0 \in \mathbb{R}^d$, and then mapped back to the conductivity space using the decoder,
\[
  \hat{K} (\bm{z}_0;\psi)=D_{\psi}(\bm{z}_0).
\]
where $\hat{K}$ represents the reconstructed approximation of the original input $K$.

\subsubsection*{VAE loss function}
The VAE is trained to reconstruct input fields while enforcing a well-structured latent space.
Assume a standard normal prior
$p(\bm{z}_0)=\mathcal{N}(\bm{0},\bm{I})$, the model is trained by minimizing the following loss:
\begin{equation}\label{eq:vae_loss}
  \mathcal{L}_{\text{VAE}}(\phi,\psi)
  =\bigl\lVert K-\hat{K} (\bm{z}_0;\psi)\bigr\rVert_{1}
  +\lambda_{\text{KL}}\,
   D_{\text{KL}}\!\bigl(q_{\phi}(\bm{z}_0\!\mid K)\,\Vert\,\mathcal{N}(\bm{0},\bm{I})\bigr),
\end{equation}

where $\phi$ and $\psi$ are the trainable parameters, respectively,
$D_{\text{KL}}(\cdot\Vert\cdot)$ denotes the Kullback-Leibler divergence, and $\lambda_{\text{KL}}$ 
is a weighting factor that balances reconstruction fidelity and latent regularization.
The first term encourages accurate recovery of the input data, while the KL divergence regularizes the latent posterior to remain close to the prior, making the learned space smooth, compact, and suitable for sampling in the subsequent diffusion model.
Specifically, the $\ell_1$-norm is used for reconstruction loss to promote sparsity and preserve sharp contrasts in the reconstructed conductivity fields, which is particularly important for capturing geological features such as sharp interfaces and channelized structures.

After training, the encoder defines the latent manifold in which the diffusion process (formulated in Eqs.~\eqref{eq:forward_transition} and \eqref{eq:ddim_update}) is applied.

%%%%%%%%%%%%%%%%%%%%%%%%%%%%

\subsubsection{LDM training and sampling}\label{sec:ldm_sample}
\begin{figure}[htb]
  \centering
  \begin{subfigure}[b]{0.4\textwidth}
    \centering
    \includegraphics[width=\textwidth]{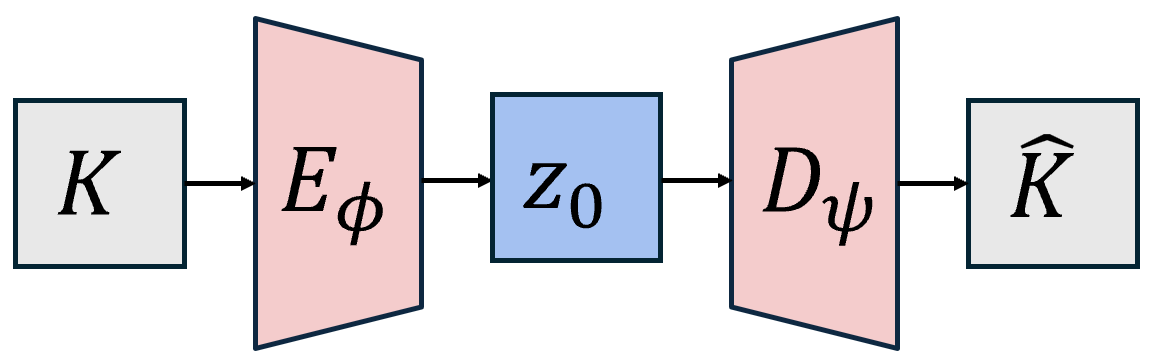}
    \caption{Variational autoencoder (VAE)}
    \label{fig:vae}
  \end{subfigure}
  \vspace{1em}
  \begin{subfigure}[b]{0.6\textwidth}
    \centering
    \includegraphics[width=\textwidth]{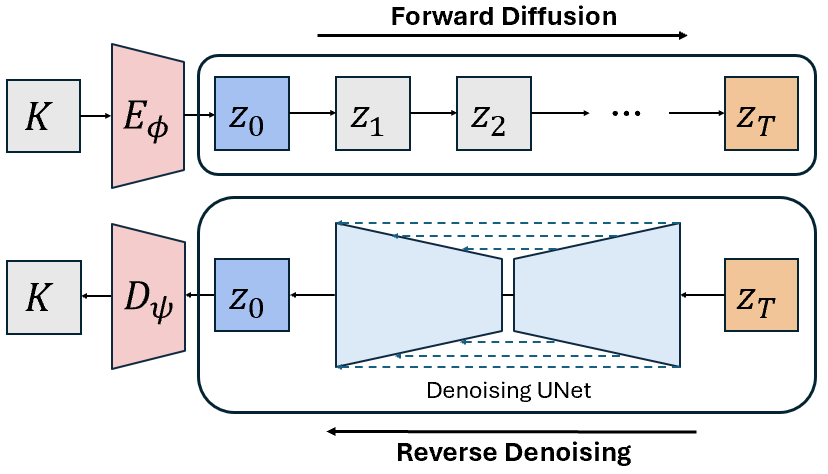}
    \caption{Diffusion process in latent space}
    \label{fig:diffusion_process}
  \end{subfigure}
  \caption{Training workflow of the latent diffusion model (LDM). The LDM comprises two key components: (a) a variational autoencoder (VAE) that compresses conductivity fields into low-dimensional latent representations, and 
  (b) a diffusion process performed in the latent space, consisting of a forward diffusion that progressively adds noise and a reverse denoising process that removes noise and reconstructs clean latent representations using a denoising U-Net. 
  Training is conducted in two stages: the VAE is first trained and fixed, followed by training the U-Net model to learn the reverse denoising dynamics in the latent space.}
  \label{fig:ldm_training}
\end{figure}

The training procedures to construct a LDM surrogate are carried out in two sequential stages, as shown in Fig.~\ref{fig:diffusion_process}.
First, 
given a dataset of conductivity field samples,
we  train a VAE using the loss function defined in \eqref{eq:vae_loss} to 
learn a low-dimensional latent representation
(Fig.~\ref{fig:vae}). 
The encoder maps each conductivity field $K$ to a latent code $z_0 = E_{\phi}(K)$ following Eq.~\eqref{eq:vae_posterior}, while the decoder reconstructs the field as $\tilde{K} = D_{\psi}(z_0)$. 
Second, once the VAE (i.e., $E_{\phi}$ and $D_{\psi}$) is trained and fixed,
a diffusion model is trained in the learned latent space (Fig.~\ref{fig:diffusion_process}). 
As described in Section~\ref{sec:diff_model}, a U-Net denoiser $\varepsilon_\theta(z_t,t)$
is trained to model the reverse diffusion process by minimizing the noise-prediction loss defined in Eq.~\eqref{eq:diffusion_loss}.

During sampling (i.e., the generation phase), the model initializes with a random latent vector  $z_T \sim \mathcal{N}(\mathbf{0}, \bm{I})$ 
and iteratively denoises it using the DDIM update (Eq.~\eqref{eq:ddim_update}) to obtain a clean latent code $z_0$.
This latent code is then decoded by
$D_{\psi}(z_0)$ to generate a conductivity realization $\hat{K}$.

% ---------------------------------------------------------------
% FVM
% ---------------------------------------------------------------

\subsection{Differentiable finite volume solver}
\label{sec:fvm}

To solve the steady-state Darcy flow equation in Eq.~\eqref{eq:darcy_governing} given a reconstructed conductivity field, we employ a finite volume method (FVM)~\citep{hirsch2007numerical, eymard2000finite}.
The core idea of  FVM is to integrate the governing equation over a finite control volume $\mathcal{V}_l$, corresponding a grid cell:
\begin{equation}
\int_{\mathcal{V}_l} \nabla \cdot (K \nabla h) dV = 0.
\end{equation}
Applying the divergence theorem, the volume integral of the divergence term is converted into a surface integral over the cell boundary $\partial\mathcal{V}_l$:

\begin{equation}\label{eq:FVM_integral}
    \int_{\mathcal{V}_l} \nabla \cdot
    (
    K
    \nabla h) dV = \oint_{\partial\mathcal{V}_l} (K\nabla h) \cdot \bm{n} \, dS = 0
\end{equation}
where $\bm{n}$ is the outward unit normal vector.
Note that the face fluxes are approximated using a two-point flux approximation (TPFA), 
which models each flux as the product of a face transmissibility and the head jump.
This cell-centered TPFA scheme preserves local conservation and achieves first-order accuracy on Cartesian grids~\citep{eymard2000finite, versteeg2007introduction}.

Discretizing the domain using a uniform Cartesian grid with $N$ degrees of freedom, the FVM yields a sparse linear system:
\begin{equation}\label{eq:fvm_solv}
\bm{A}(\hat{\bm K})\,\hat{\bm{h}}=\bm{b}.
\end{equation}
where $\hat{\bm{K}}, \hat{\bm{h}} \in \mathbb{R}^N$ are the discrete conductivity and pressure head fields, $\bm{A} \in \mathbb{R}^{N\times N}$ 
is the stiffness matrix, and $\bm{b} \in \mathbb{R}^N$ encodes the boundary conditions. 
Dirichlet boundary conditions are imposed by replacing the corresponding matrix rows with identity rows and assigning the prescribed head values to the right-hand side~\citep{patankar2018numerical, versteeg2007introduction},
while
homogeneous Neumann (no-flux) boundaries are handled by zeroing out the normal transmissibility on boundary faces~\citep{eymard2000finite}. 
A detailed derivation of the 2D FVM discretization for heterogeneous isotropic media is provided in ~\ref{sec:app_A}.

\subsubsection*{Differentiable solver}
To support gradient-based inversion, the FVM solver is implemented in 
\texttt{JAX}~\citep{frostig2019compiling, jax2018github}, 
enabling efficient computation of reverse-mode gradients via automatic differentiation.
This makes the discrete forward operator
$\hat{\mathcal{M}}: \hat{\bm K} \mapsto \hat{\bm h}$ fully differentiable,
allowing for exact  adjoint-based sensitivity analysis. 
Compared to conventional finite-difference sensitivity analysis, 
this differentiable solver 
ensures seamless integration with the LDM framework
and supports scalable optimization of latent variables 
within the coupled PDE-constrained system.

% ---------------------------------------------------------------
% LD-DIM
% ---------------------------------------------------------------
\subsection{LD-DIM: Physics-constrained inverse optimization}\label{sec:physics_inverse}
\begin{figure}[htb!]
  \centering
  \includegraphics[width=0.85\textwidth]{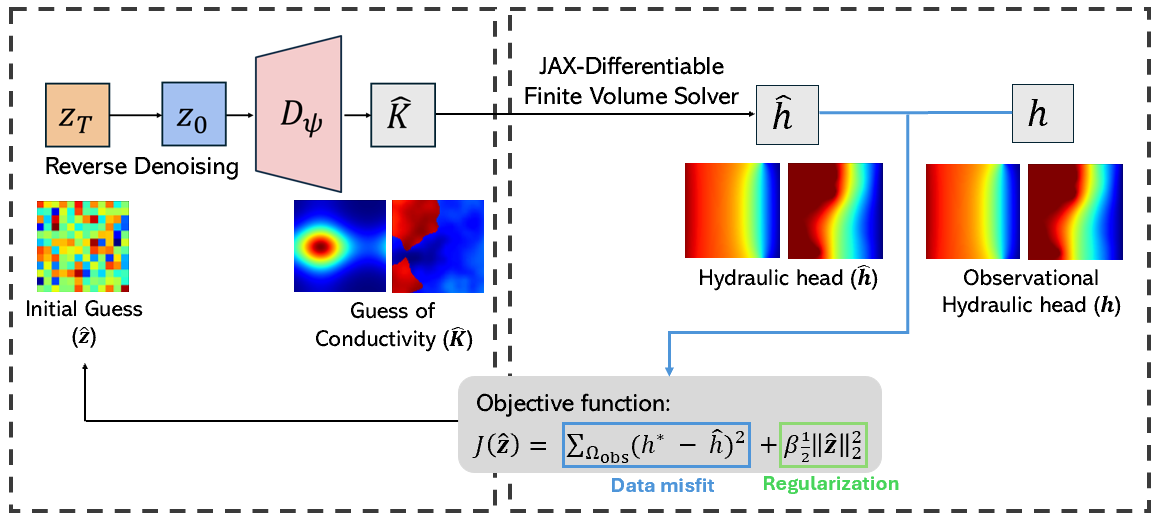}
  \caption{LD-DIM inverse modeling workflow. The pretrained latent diffusion model maps a low-dimensional latent vector $\hat{\bm{z}}$ into a conductivity field via reverse denoising and decoding. The generated conductivity field is input to a differentiable finite-volume solver that computes the hydraulic head $\hat{h}$, which is compared with observations $h^{*}$ in the objective function. Gradients with respect to latent variable $\hat{\bm{z}}$ are efficiently computed by combining adjoint-based gradients through the finite-volume solver and automatic differentiation through the pre-trained latent diffusion model, enabling end-to-end optimization in the compact latent space.}

  \label{fig:inverse_framework}
\end{figure}

The proposed LD-DIM framework integrates a pretrained LDM (Section \ref{sec:ldm}) with the differentiable finite-volume solver (Section \ref{sec:fvm}) to reconstruct spatially varying parameter fields from sparse observations (Fig. \ref{fig:inverse_framework}). 
By treating the LDM as a learned prior, the resulting inversion is formulated as a deterministic optimization problem in a low-dimensional latent space. This formulation constrains the inverse solution to a physically admissible latent manifold while rigorously enforcing the governing equations through the embedded physics-based solver.
% this approach constrains the inverse problem within a geologically plausible latent manifold while ensuring strict satisfaction of governing equations through an embedded physics-based  solver.
Furthermore, an efficient optimization procedure (Algorithm~\ref{alg:pipeline})  is developed that operates entirely in the low-dimensional latent space, combining automatic differentiation through the generative model with adjoint-based gradient computation for fully end-to-end optimization.

Given an initial latent guess $\hat{\bm{z}} \sim \mathcal{N}(\bm{0},\bm{I})$, the LDM generates a conductivity field $\hat{\bm{K}} = \mathcal{G}(\hat{\bm{z}}; \theta, \psi)$, where $\mathcal{G}$ denotes the generative 
sampling operator involving the learned reverse diffusion process $\varepsilon_\theta$ and the trained decoder $D_{\psi}$.
The details of the generative process are provided in Section~\ref{sec:ldm_sample}.

With the pretrained LDM surrogate $\hat{\bm{K}} = \mathcal{G}(\hat{\bm{z}})$, the LD-DIM framework is formulated
as the following  PDE-constrained optimization problem:
\begin{align}\label{eq:opt_DIM}
\begin{split}
& \min_{\hat{\bm{z}}} \ \mathcal{J}(\hat{\bm{z}},\hat{\bm{K}}(\hat{\bm{z}}),\hat{\bm{h}}) \\
& \text{s.t.} \quad 
  \bm{A}(\hat{\bm{K}}(\hat{\bm{z}}))\,\hat{\bm{h}}= \bm{b}
\end{split}
\end{align}
Here, the constraint is enforced using the differentiable finite-volume solver introduced earlier, and the objective function $\mathcal{J}$ is detailed in the following subsection.

%%%%%%%%%%%%%%%%%%%%%%%%%%%%%%
\subsubsection{Pretrained stage}

As outlined in Section~\ref{sec:ldm}, the LDM is pretrained on synthetic conductivity fields. The VAE is first trained using Eq.~\eqref{eq:vae_loss} to learn the encoder $E_{\phi}$ and decoder $D_{\psi}$, followed by training the diffusion model in the latent space using Eq.~\eqref{eq:diffusion_loss}. 

During inversion, the pretrained LDM defines a generative mapping
$\mathcal{G}(\cdot)$ that transforms latent variables  $\hat{\bm{z}}$ into corresponding conductivity fields $\hat{\bm{K}}$. 
Therefore, the optimization is carried out with respect to the latent space, ensuring that the inferred parameters remain confined to the geologically consistent manifold learned from the synthetic training data.

%%%%%%%%%%%%%%%%%%%%%%%%%%%%%%
\subsubsection{Objective function}
The inversion minimizes the objective function:
\begin{equation}\label{eq:tot_loss}
  \mathcal{J}(\hat{\bm{z}}) = \ell(\hat{\bm{z}}) + \beta\,\mathcal{R}(\hat{\bm{z}})
\end{equation}
where $\ell$ measures the data misfit and $\mathcal{R}$ is a regularization term.
The misfit term is defined as:
\begin{equation}\label{eq:data_misfit}
  \ell(\hat{\bm{z}}) = \sum_{i \in \mathcal{I}_{\text{obs}}} \bigl(h^{*}_i-\hat{h}_i\bigr)^{2}
\end{equation}
with $h^{*}_i$ and $\hat{h}_i$
denoting the
measured and predicted hydraulic heads at observation locations $\{\bm{x}_i\}_{i \in \mathcal{I}_{\text{obs}}}$, respectively.
The regularization term
\begin{equation}\label{eq:reg}
  \mathcal{R}(\hat{\bm{z}}) = \frac{1}{2} \lVert\hat{\bm{z}}\rVert_2^{2}
\end{equation}
encourages the latent variables to remain close to the standard normal prior, with $\beta$
controlling the regularization strength.

%%%%%%%%%%%%%%%%%%%%%%%%%%%%%%

\subsubsection{Optimization procedure \& sensitivity analysis}

Leveraging the end-to-end differentiable framework, we employ efficient hybrid gradient-based optimization to solve for the optimal latent variables in 
Eq. \eqref{eq:opt_DIM}.
The sensitivity of the objective function $\mathcal{J}(\hat{\bm{z}})$ in Eq.~\eqref{eq:tot_loss} with respect to the latent variable $\hat{\bm{z}}$ is computed via the chain rule:
\begin{equation}\label{eq:chain_rule}
  \nabla_{\hat{\bm{z}}}\mathcal{J} = \left(\frac{\partial\hat{\bm K}}{\partial\hat{\bm{z}}}\right)^{\!\top} \nabla_{\hat{K}}\ell + \beta\hat{\bm{z}}
\end{equation}
where $\partial \hat{\bm K}/\partial \hat{\bm{z}}$ is computed 
through \textit{automatic differentiation} with the pretrained LDM.

Since the simulated pressure field $\hat{\bm{h}}$ depends on the conductivity field $\hat{\bm{K}}$ through the PDE constraint, the gradient $\nabla_{\hat{\bm{K}}} \ell$ is given by
\begin{equation}\label{eq:gradK_l}
  \nabla_{\hat{K}}\ell = \left(\frac{\partial\hat{\bm h}}{\partial\hat{\bm{K}}}\right)^{\!\top} 
  \nabla_{\hat{h}}\ell
\end{equation}
To compute this efficiently, we adopt the \textit{discrete adjoint method}~\cite{dwight2006effect,plessix2006review}
%giles2000introduction
to compute this sensitivity efficiently. Letting $\bm{r} = \nabla_{\hat{\bm{h}}} \ell$, the adjoint equation is
\begin{equation}\label{eq:adjoint}
  \bm{A}(\hat{\bm K})^{\!\top}\bm{\lambda}= \bm{r}
\end{equation}
where $\bm{\lambda} \in \mathbb{R}^N$ is the discrete adjoint state. Once $\bm{\lambda}$ is obtained, the conductivity gradient is computed as
\begin{equation}\label{eq:grad_k}
  \nabla_{\hat{K}}\ell
  = -\,
    \bm{\lambda}^{\!\top}
    \Bigl(\tfrac{\partial\bm{A}(\hat{\bm{K}})}{\partial\hat{\bm{K}}}\Bigr)
    \hat{\bm{h}}
\end{equation}
where $\frac{\partial \bm{A}(\hat{\bm{K}})}{\partial \hat{\bm{K}}}$ represents the sensitivity of the stiffness matrix with respect to the conductivity field, 
which can be constructed analytically using the local dependence of $\hat{\bm{K}}$ on each grid cell. Implementation details can be referred to~\ref{app:implementation}.

We remark that the matrix $\bm{A}$ is reused in both the forward solve and adjoint computation, %ensuring consistency and computational efficiency. 
ensuring consistency between the primal and adjoint problems while improving computational efficiency. 
This adjoint-based approach avoids the explicit construction or storage of the dense Jacobian $\frac{\partial \hat{\bm{h}}}{\partial \hat{\bm{K}}}$, thereby substantially reducing memory requirements and computational overhead.
The overall computational workflow of the proposed LD-DIM framework for estimating unknown conductivity fields is summarized in Algorithm~\ref{alg:pipeline}.

While the present study focuses on a steady-state elliptic PDE, the proposed hybrid adjoint–automatic differentiation strategy, combining discrete adjoint methods for the numerical solver with automatic differentiation through the latent generative model, is not specific to Darcy flow. The same computational paradigm applies to a broad class of PDE-constrained inverse problems with differentiable discretizations. Extensions to more complex forward models are conceptually straightforward but require problem-specific adjoint formulations, such as those reflected in Eqs.~\eqref{eq:adjoint} and~\eqref{eq:grad_k}.

\begin{algorithm}[H]
\caption{Workflow of Latent Diffusion Differentiable Inversion Method (LD-DIM)}\label{alg:pipeline}
\begin{algorithmic}[1]
\State Initialize latent vector $\hat{\bm{z}}^{(0)}\sim\mathcal{N}(\bm{0},\bm{I})$
\For{$\nu=0$ \textbf{to} maxIter}
  \State Generate conductivity field: $\hat{\bm K}\gets \mathcal{G}(\hat{\bm{z}}^{(\nu)})$
  \State Solve PDE system \eqref{eq:fvm_solv} for $\hat{\bm{h}}$
  \State Compute loss $\mathcal{J}$ via \eqref{eq:tot_loss}
  \State Compute gradient $\nabla_{\hat{\bm{z}}}\ell$ via chain rule and adjoint (Eqs. \eqref{eq:chain_rule}--\eqref{eq:grad_k})
  \State Compute total gradient: $\nabla_{\hat{\bm{z}}} \mathcal{J} = \nabla_{\hat{\bm{z}}} \ell + \nabla_{\hat{\bm{z}}} \mathcal{R}$
  \State Update latent vector: $\hat{\bm{z}}^{(\nu+1)}\gets\hat{\bm{z}}^{(\nu)}-\eta\,\nabla_{\hat{\bm{z}}}\mathcal{J}$
\EndFor
\State \textbf{return} Final conductivity estimate: $\hat{\bm K}\gets \mathcal{G}(\hat{\bm{z}}^{*})$
\end{algorithmic}
\end{algorithm}

\section{Numerical results}\label{sec:numerical-validation}

We evaluate the effectiveness of the proposed LD-DIM framework through a series of numerical experiments on PDE-constrained inverse problems for flow in porous media. 
The objective is to reconstruct spatially heterogeneous conductivity fields from sparse hydraulic head observations using two representative test cases: (1) Gaussian conductivity fields with varying correlation lengths (Section \ref{subsec:gaussian-fields}), and (2) complex bimaterial distributions with sharp interfaces (Section \ref{subsec:bimaterial-fields}). 

For each case, the LDM is pretrained on the corresponding synthetic datasets generated using the procedures described below. This setup enables a systematic  evaluation of reconstruction performance under controlled conditions.

%%%%%%%%%%%%%%%%%%%%%
%%
\subsection{Computational setup}\label{subsec:governing-physics}

All numerical experiments are based on steady-state groundwater flow governed by Darcy's equation~\eqref{eq:darcy_governing} over a 
2D computational domain $\Omega = [-50,50]  \times [-50,50] ~\text{m}$,
 associated with the Dirichlet and Neumann boundary conditions specified as:
\begin{subequations}\label{eq:boundary-conditions}
\begin{align}
  h(-50,y) = 1,\quad h(50,y) = 0, 
  \qquad &\text{for } y \in [-50,50]~\mathrm{m}\, \label{eq:bc-dirichlet}\\[2pt]
  \frac{\partial h}{\partial n}(x,-50) = 0,\quad \frac{\partial h}{\partial n}(x,50) = 0,
  \qquad &\text{for } x \in [-50,50]~\mathrm{m}
  \label{eq:bc-neumann}
\end{align}
\end{subequations}
where $\partial/\partial n$ denotes the outward normal derivative on the boundary.

The domain is discretized using a uniform $100 \times 100$ Cartesian grid (i.e., $N=10,000$), with grid spacing $\Delta x = \Delta y = 1$~m. Synthetic head observations are sampled at interior points according to uniform spatial patterns of varying density to assess the robustness of LD-DIM under different levels of data sparsity.

To assess the accuracy of reconstructed conductivity and pressure head fields, we employ 3 complementary metrics that quantify both global error magnitude and structural fidelity:
\begin{itemize}
\item \textbf{Relative $L^2$ error:} 
\begin{equation}
  %\mathcal{E}_{2} 
  \epsilon_{K}= \frac{\bigl\lVert \hat{K}-K^{*}\bigr\rVert_{2}}{\bigl\lVert K^{*}\bigr\rVert_{2}}, \qquad 
  \epsilon_{h}= \frac{\bigl\lVert \hat{h}-h^{*}\bigr\rVert_{2}}{\bigl\lVert h^{*}\bigr\rVert_{2}},
  \label{eq:relative-error}
\end{equation}
where the symbols $\widehat{(\cdot)}$ and $(\cdot)^{*}$ denote the reconstructed and reference (ground-truth) quantities, respectively, and  $\lVert\cdot\rVert_{2}$ denotes the $L^2$ norm
%applied to vectorized fields.
over the discretized grid.

\item \textbf{Mean-corrected relative $L^2$ error (for Bimaterial fields):} For two-phase bimaterial cases characterized by discontinuous interfaces, we adopt a mean-adjusted error to better capture spatial pattern mismatch:
\begin{equation}
  \tilde{\epsilon}_{K} =
    \frac{\lVert\hat{K}-K^{*} \rVert_{2}}
         {\lVert K^{*} -\bar K^{*}\rVert_{2}},  
  \label{eq:relative-error-bimaterial}
\end{equation}
where $\bar K^{*}$ is the spatial mean of $K^{*}$. 
This error metric is used exclusively for bimaterial
cases in Section~\ref{subsec:bimaterial-fields}.

\item \textbf{Structural Similarity Index (SSIM):} Originally developed for image quality assessment~\cite{wang2004image}, the SSIM metric measures the structural similarity between two images by considering luminance, contrast, and structure. For arbitrary two random fields, $p$ and $q$,  the SSIM metric is defined as:
\begin{equation}
  \epsilon_{\text{SSIM}}(p,q) = 
\frac{\bigl(2\mu_{p}\mu_{q}+C_{1}\bigr)\bigl(2\sigma_{pq}+C_{2}\bigr)}{\bigl(\mu_{p}^{2}+\mu_{q}^{2}+C_{1}\bigr)\bigl(\sigma_{p}^{2}+\sigma_{q}^{2}+C_{2}\bigr)},
  \label{eq:ssim}
\end{equation}
where $\mu_p$ and $\mu_q$ are local means, $\sigma_p^2$ and $\sigma_q^2$ are local variances, and $\sigma_{pq}$ is the local covariance between $p$ and $q$. 
In our context, $p$ and $q$ correspond to the ground-truth and reconstructed conductivity fields, respectively.
The stabilizing constants are $C_1 = (k_1 L)^2$ and $C_2 = (k_2 L)^2$, with $L$ denoting the dynamic range and default values $k_1 = 0.01$, $k_2 = 0.03$. SSIM values range from $-1$ to $1$, with higher scores indicating greater structural similarity. This metric is especially useful for evaluating the spatial coherence of reconstructed conductivity fields~\citep{brunet2011mathematical}. 
\end{itemize}

Unless otherwise specified, all error metrics are evaluated over the entire grid.

\subsection{Case study 1: Gaussian conductivity fields}\label{subsec:gaussian-fields}
We begin with a benchmark scenario involving Gaussian random fields to represent heterogeneous conductivity distributions.
Section~\ref{subsec:gaussian-generation} outlines the synthetic data generation process and the training setup for the latent diffusion model.
Next, we present quantitative and visual inversion results (Section~\ref{subsec:gaussian-results}) 
and compare LD-DIM with a physics-informed neural network (PINN) baseline that directly enforces Darcy’s law during training (Section~\ref{subsec:pinn-comparison}). 
Through these analyses, we assess the reconstruction accuracy, structural consistency, and robustness of LD-DIM across different correlation lengths and observation densities.

%%%%%%%%%%%%%%
\subsubsection{Synthetic data generation \& LDM Training}\label{subsec:gaussian-generation}
 
In this example, the spatially correlated conductivity fields are generated using a spectral synthesis procedure~\cite{shinozuka1996simulation}.
The log-conductivity field $Y(\bm{x})$ is constructed through applying an inverse Fourier transformation to spectrally filtered white noise:
\begin{equation}
\begin{aligned}
& K(\bm{x}) = \exp\!\bigl[\,Y(\bm{x})\,\bigr], \\
& \text{where} \quad Y(\bm{x}) = 
\mathcal{F}^{-1}\!\left\{\sqrt{S(\bm{k})}\,\mathcal{F}\{\xi\}\right\}, \quad
S(\bm{k}) &= 
\exp\!\Bigl(-\tfrac{1}{2}\lVert\bm{k}\rVert_{2}^{2}\,\lambda^{2}\Bigr),
\end{aligned}\label{eq:K-definition}
\end{equation}
Here, $\xi \sim \mathcal{N}(0,1)$ is complex white noise, $\mathcal{F}$ and $\mathcal{F}^{-1}$ denote discrete Fourier and inverse Fourier transforms,
$\bm{k}$ is the wavenumber vector,
and $\lambda$ controls the spatial correlation length. 
Increasing $\lambda$  leads to smoother and more spatially coherent fields by suppressing high-frequency components.
To ensure numerical stability, the log-field \(Y\) in Eq.~\eqref{eq:K-definition} is standardized to zero mean and unit variance, then clipped to the interval \([-3,3]\) before applying the exponential transformation, the resulting conductivity field $K$ is subsequently normalized to the range $[0,1]$.

\begin{figure}[htb!]
  \centering
  \includegraphics[width=0.6\textwidth]{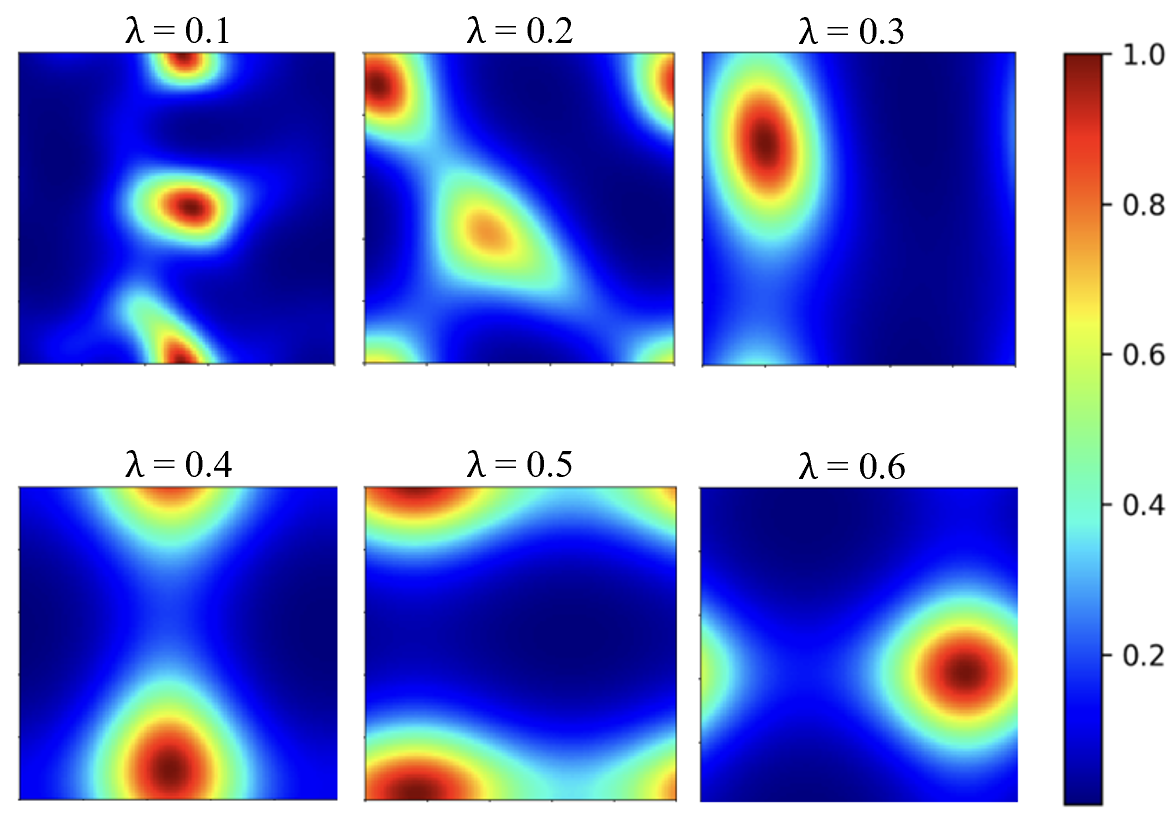}
  \caption{Representative samples of Gaussian conductivity fields corresponding to different correlation lengths $\lambda$.
  Top row (left to right): $\lambda=0.1,\,0.2,\,0.3$. 
  Bottom row (left to right): $\lambda=0.4,\,0.5,\,0.6$. 
Fields exhibit progressively smoother spatial structure with increased correlation length.}
  \label{fig:gaussian-correlation-samples}
\end{figure}

A total of 2,000 conductivity field samples are generated using Eq.~\eqref{eq:K-definition}, evenly distributed across seven correlation lengths $\lambda \in [0.1, 0.7]$ to capture a broad spectrum of spatial variability. Of these, 1,400 samples are randomly selected to train the LDM, 400 for LDM validation, and 
the remaining 200 are reserved for inverse modeling tests.
Figure~\ref{fig:gaussian-correlation-samples} shows representative samples, highlighting the transition from highly heterogeneous fields to smoother spatial structures as $\lambda$ increases.

The LDM prior is trained on the generated fields using the approach described in Section~\ref{sec:ldm_sample}. 
A variational autoencoder compresses each conductivity field into a $12 \times 12$ latent representation, followed by diffusion training in the latent space. This enables efficient sampling and supports gradient-based inversion within a geologically informed manifold. 

\begin{figure}[htbp!]
  \centering
  \includegraphics[width=0.95\textwidth]{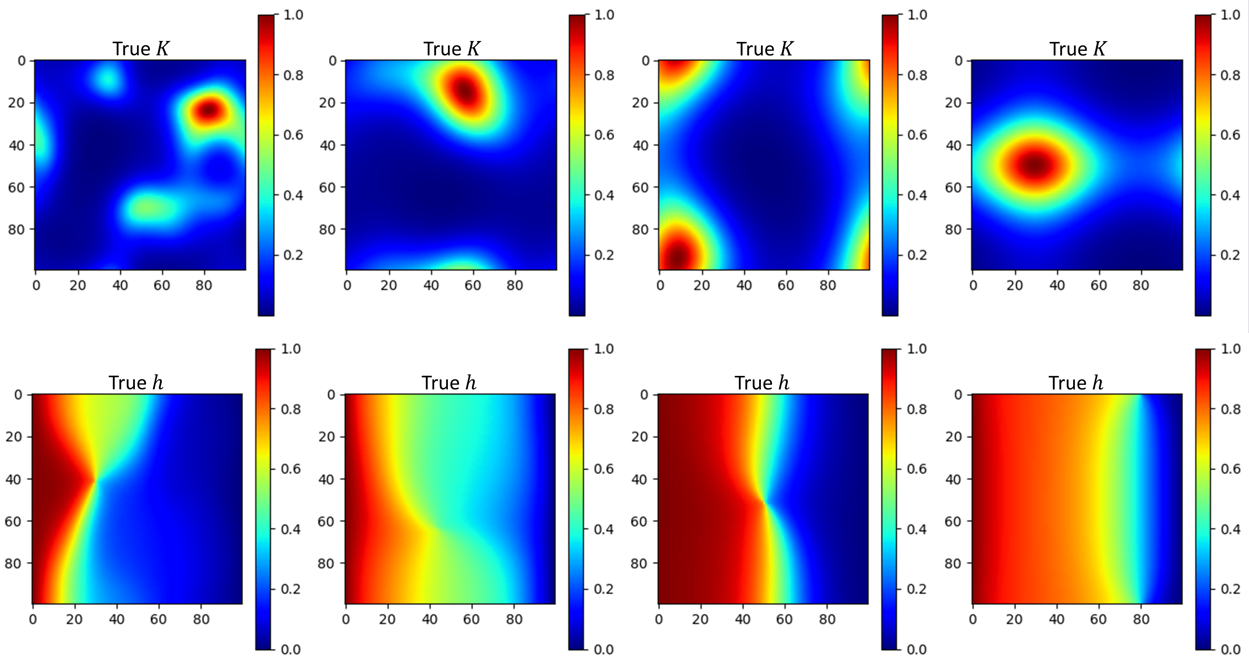}
  \caption{Representative testing samples of Gaussian conductivity fields $K$ (top) and the corresponding hydraulic head fields $h$ (bottom) for correlation lengths $\lambda = 0.1,\,0.2,\,0.3,$ and $0.4$ (left to right). 
  These test data serve as the ground truth for the inverse experiments.
  }
  \label{fig:Gaussian_reference_K_h}
\end{figure}

For the following inversion tests, reference hydraulic head solutions are computed using the FVM solver introduced in Section~\ref{sec:fvm}, based on the corresponding conductivity fields.
Some ground-truth conductivity–head pairs used in inverse tests for $\lambda = 0.1$, $0.2$, $0.3$, and $0.4$ are displayed in Fig. ~\ref{fig:Gaussian_reference_K_h}.

%%%%%%%%%%%%%%%%%%%%%%%%%%%
%%
\subsubsection{Conductivity characterization results}\label{subsec:gaussian-results}

We evaluate the ability of LD-DIM to recover various testing conductivity fields from sparse hydraulic head observations using the reference solutions illustrated in Fig.~\ref{fig:Gaussian_reference_K_h}, where 256 hydraulic head observations arranged in a $16 \times 16$ uniform grid are used. 
We note that these observations are generated from conductivity fields not seen during training, thereby assessing out-of-distribution inversion capability.

\begin{figure}[htbp!]
  \centering

  % --- (a) λ = 0.1 ---
  \begin{subfigure}[b]{0.8\textwidth}
    \centering
    \includegraphics[width=\textwidth]{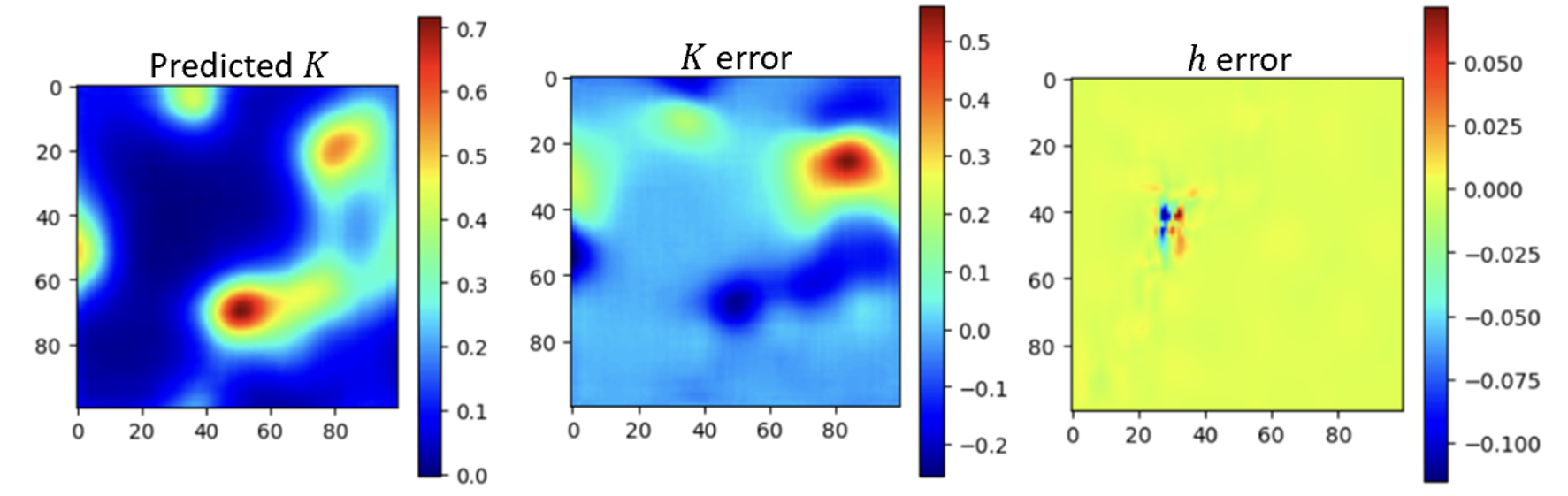}
    \caption{Correlation length $\lambda = 0.1$. $\epsilon_{K}$: $4.87 \times 10^{-1}$; SSIM: 0.709.}
    \label{fig:gaussian-lambda-01}
  \end{subfigure}

  % --- (b) λ = 0.2 ---
  \begin{subfigure}[b]{0.8\textwidth}
    \centering
    \includegraphics[width=\textwidth]{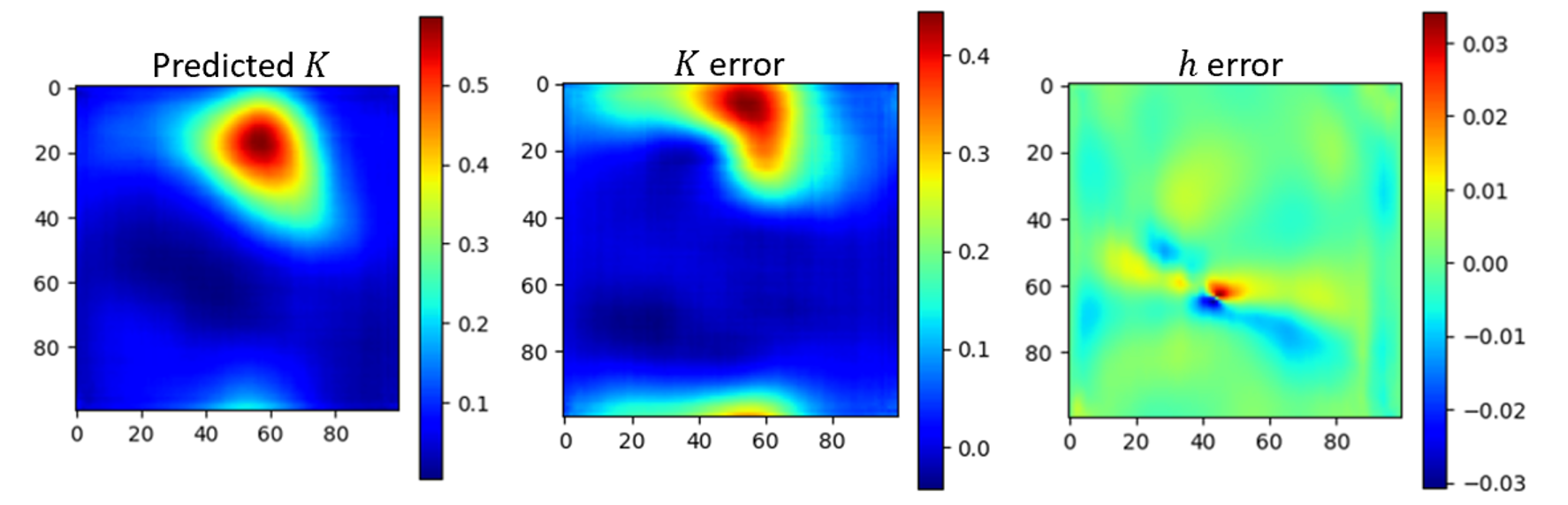}
    \caption{Correlation length $\lambda = 0.2$. $\epsilon_{K}$: $4.031 \times 10^{-1}$; SSIM: 0.760.}
    \label{fig:gaussian-lambda-02}
  \end{subfigure}

  % --- (c) λ = 0.3 ---
  \begin{subfigure}[b]{0.8\textwidth}
    \centering
    \includegraphics[width=\textwidth]{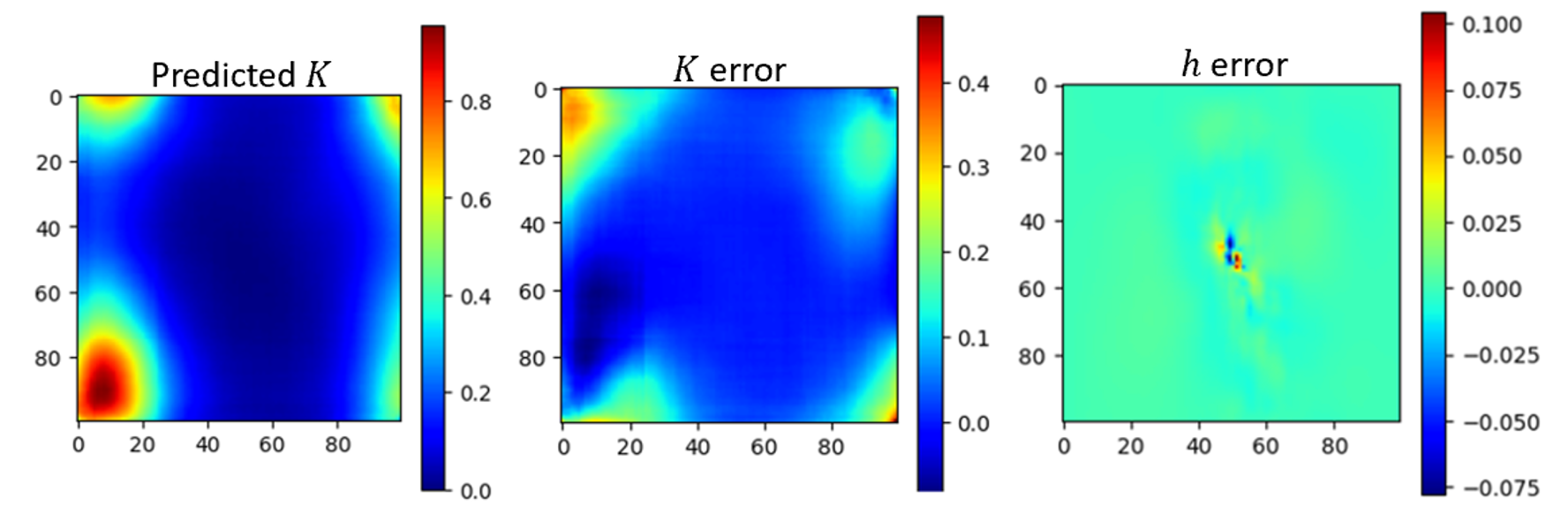}
    \caption{Correlation length $\lambda = 0.3$. $\epsilon_{K}$: $2.642 \times 10^{-1}$; SSIM: 0.917.}
    \label{fig:gaussian-lambda-03}
  \end{subfigure}

  % --- (d) λ = 0.4 ---
  \begin{subfigure}[b]{0.8\textwidth}
    \centering
    \includegraphics[width=\textwidth]{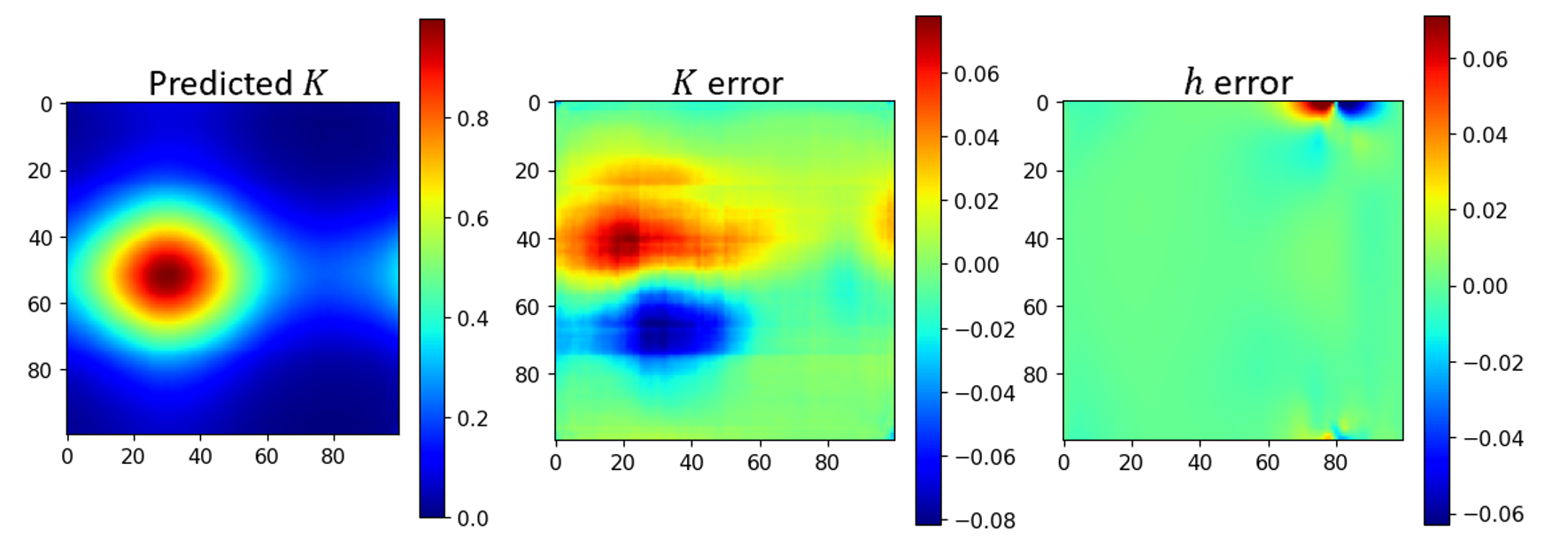}
    \caption{Correlation length $\lambda = 0.4$. $\epsilon_{K}$: $3.193 \times 10^{-2}$; SSIM: 0.987.}
    \label{fig:gaussian-lambda-04}
  \end{subfigure}

  \caption{Inversion results by LD-DIM for Gaussian conductivity fields with varying correlation lengths. Each row displays (left to right): predicted conductivity field, absolute error in conductivity, and error in hydraulic head. The framework demonstrates improved accuracy for larger correlation lengths due to reduced high-frequency content and enhanced information content in sparse observations.}
  \label{fig:gaussian-inversion-results}
\end{figure}

Fig.~\ref{fig:gaussian-inversion-results} presents representative inversion results across correlation lengths from $\lambda = 0.1$ to $0.4$, with each row showing the reconstructed conductivity field, and the corresponding absolute error in $K$ and hydraulic head $h$.
With a fixed number of $h$ observations, the quality of inversion is strongly influenced by the spatial correlation length. 
Nevertheless, LD-DIM successfully captures the overall structure of the underlying conductivity fields across all cases. 
The conductivity reconstruction error $\epsilon_K$ decreases significantly from $4.87 \times 10^{-1}$ at $\lambda = 0.1$ to $3.19 \times 10^{-2}$ at $\lambda = 0.4$, while SSIM scores increase from 0.709 to 0.987.
These trends highlight the influence of spatial correlation: larger $\lambda$ leads to more coherent patterns that are easier to infer, whereas small $\lambda$ presents significant challenges due to limited data support for fine-scale heterogeneity.

In addition to conductivity, LD-DIM achieves high-quality predictions of the full hydraulic head field. As shown in Fig.~\ref{fig:gaussian-inversion-results}, the pointwise error in $h$ remains below 10\% across all cases, regardless of correlation length. The largest discrepancies in $h$ tend to occur in regions with highly localized head variations, especially in the $\lambda = 0.1$ and $\lambda = 0.3$ cases (see Fig. \ref{fig:Gaussian_reference_K_h}). Despite this, the overall predictive performance remains robust, confirming the effectiveness of LD-DIM in recovering both field structure and flow response from limited head data alone.

The ability of LD-DIM to infer physically plausible patterns under such sparse supervision underscores the critical role of the pretrained LDM as a learned geological prior. 
This prior regularizes the solution space, enabling stable and meaningful inversion even in severely ill-posed settings. 
Unlike conventional PDE-constrained or PINN-based methods (e.g., Fig.~\ref{fig:ldm-vs-pinn-l01} and Fig.~\ref{fig:ldm-vs-pinn-l04}), which often in the absence of explicit conductivity data, LD-DIM avoids numerical instability or over-smoothing by constraining the solution within a geologically plausible latent manifold.

\begin{figure}[htbp!]
  \centering
  \includegraphics[width=0.6\textwidth]{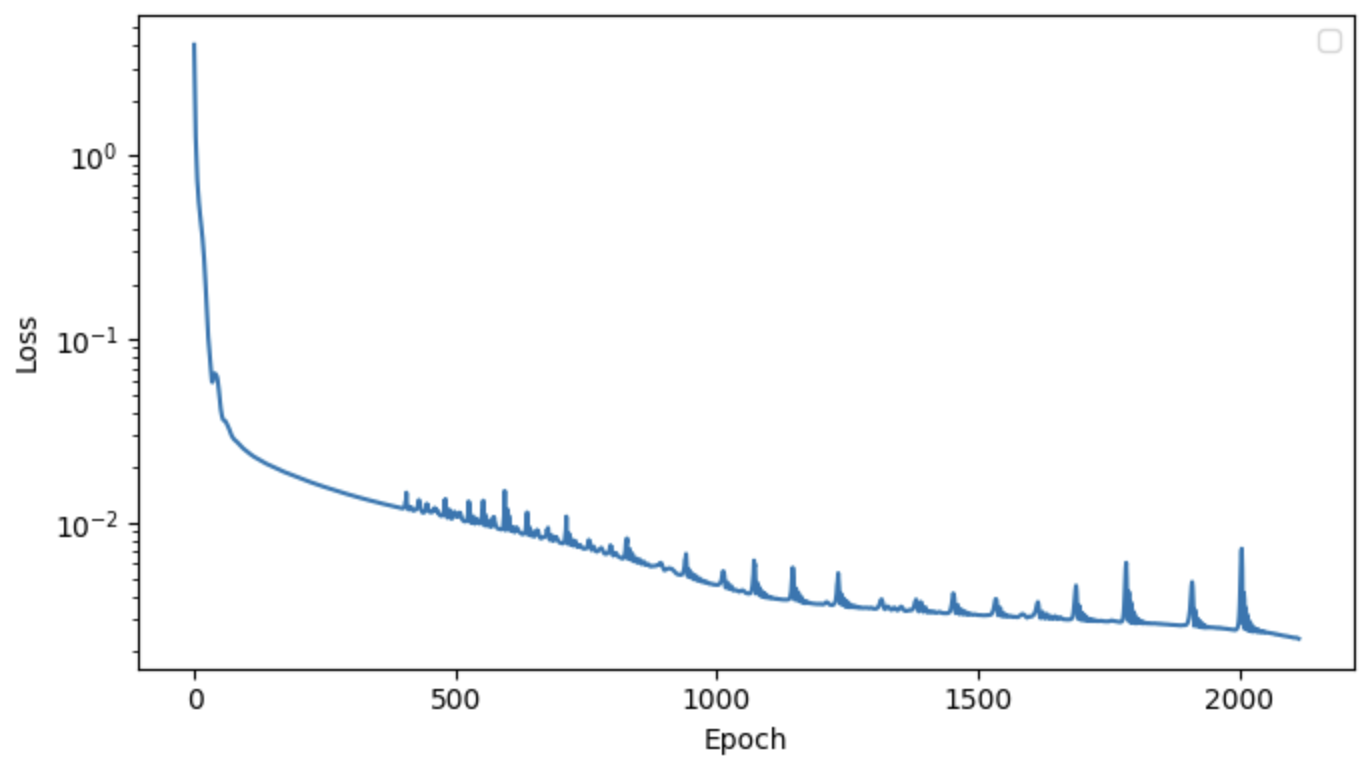}
  \caption{Convergence behavior of the inversion optimization for the case shown in Fig.~\ref{fig:gaussian-lambda-04} ($\lambda=0.4$). The data misfit loss decreases rapidly within 500 epochs and stabilizes thereafter, indicating successful convergence to a data-consistent solution.}
  \label{fig:inversion-loss}
\end{figure}

Moreover, Fig.~\ref{fig:inversion-loss} shows the loss trajectory for the representative case with $\lambda = 0.4$, illustrating the efficient and stable convergence of gradient-based optimization within the latent space under the LD-DIM framework.

In Section \ref{sec:including_k_data}, we 
further extend the framework to a hybrid setting that incorporates sparse direct conductivity observations. This extension enhances reconstruction accuracy, particularly for cases with small correlation lengths, where limited supervision would otherwise hinder performance.

%%%%%%%%%%%%%%%%%%%%%%%%%%%%%%%%%%%
%%%%%%%%%%%%%%%%%%%%%%%%%%%%%%%%%%%

%%%%%%%%%%%%%%%%%%%%%%%%%%%%%
%%
\subsubsection{Comparison between LD-DIM and physics-informed neural networks (PINNs)}\label{subsec:pinn-comparison}

Here, we introduce a PINN-based approach for comparison, which enforces Darcy’s law as a soft constraint in the loss function~\cite{raissi2019physics,he2020physics}. 
In this method, both the conductivity field and the hydraulic head are modeled using fully connected neural networks, each comprising three hidden layers with 60 neurons and hyperbolic tangent activation functions. The networks are trained by minimizing a composite loss function that combines data misfit with the residuals of the governing PDE and boundary conditions:
\begin{equation}
\label{eq:pinn-loss}
\begin{aligned}
\mathcal{L}_{\mathrm{PINN}}
&= \frac{1}{N_{\mathrm{obs}}}
   \sum_{i\in\mathcal{I}_{\mathrm{obs}}}
   \bigl(\hat h(\boldsymbol{x}_i)-h(\boldsymbol{x}_i)\bigr)^2 \\
&\quad+ w_b\,\frac{1}{N_b}
   \sum_{i\in\mathcal{I}_b}
   \bigl(\mathcal{B}[\hat h](\boldsymbol{x}_i)-b(\boldsymbol{x}_i)\bigr)^2 \\
&\quad+ w_f\,\frac{1}{N_f}
   \sum_{i\in\mathcal{I}_f}
   \Bigl[\nabla\!\cdot\!\bigl(\hat K(\boldsymbol{x}_i)\,\nabla \hat h (\boldsymbol{x}_i)\bigr)\Bigr]^2 \\
&\quad+ w_{\mathrm{reg}}\,\mathcal{R}(\Theta).
\end{aligned}
\end{equation}
where $\mathcal{I}_{\mathrm{obs}}$, $\mathcal{I}_b$, and $\mathcal{I}_f$ denote the sets of observation, boundary, and interior collocation points, with cardinalities $N_{\mathrm{obs}} = |\mathcal{X}_{\mathrm{obs}}|$, $N_b = |\mathcal{X}_b|$, and $N_f = |\mathcal{X}_f| = 1000$, respectively. The operator $\mathcal{B}[\cdot]$ represents the application of Dirichlet or Neumann boundary conditions, with $b(\bm{x})$ corresponding to the prescribed boundary values. The predictions $\hat h$ and $\hat K$ are outputs of the respective neural networks, $\Theta$ denotes the collection of network parameters, and $\mathcal{R}(\Theta)$ is a $L^2$ regularization term weighted by $w_{\mathrm{reg}}$.
The coefficients $w_b$ and $w_f$ balance the contributions of the boundary condition and PDE residual terms within the total loss.
Following \cite{he2020physics}, we set $w_b = w_f = 1$ and $w_{\mathrm{reg}} = 10^{-6}$. 
Importantly, unlike the original formulation in~\cite{he2020physics}, 
Eq.~\eqref{eq:pinn-loss} does not include direct measurements of the conductivity field $K$, aligning the experiment setup with the LD-DIM inversion in Section \ref{subsec:gaussian-results}.

% --- λ = 0.1 ---
\begin{figure}[htbp!]
  \centering
  \includegraphics[width=0.85\textwidth]{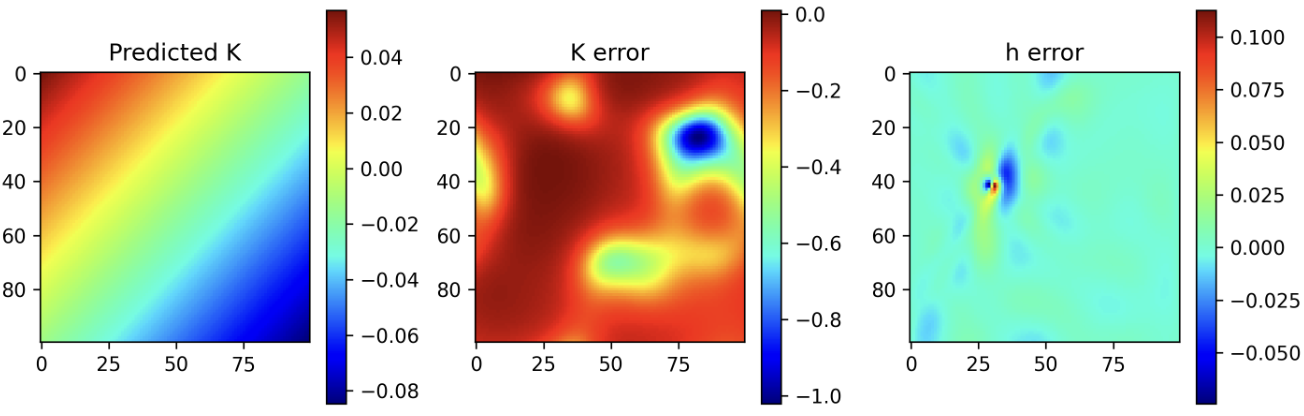}
  \caption{Inverse modeling results from PINN for a conductivity field with correlation length $\lambda=0.1$. 
  $\epsilon_{K}=1.064$, SSIM = $- 0.054$. 
  From left to right: predicted conductivity field, conductivity error, and hydraulic head error.}
  \label{fig:ldm-vs-pinn-l01}
\end{figure}

% --- λ = 0.4 ---
\begin{figure}[htbp!]
  \centering
  \includegraphics[width=0.85\textwidth]{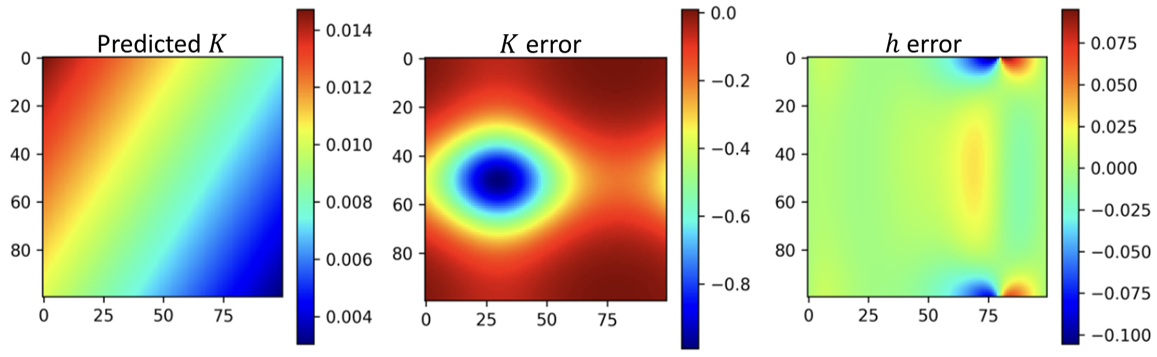}
  \caption{Inverse modeling result by PINN at correlation length $\lambda=0.4$.
  $\epsilon_{K}=9.793\times10^{-1}$, SSIM = 0.208. 
  From left to right: predicted conductivity field, error in conductivity, and error in hydraulic head.}
  \label{fig:ldm-vs-pinn-l04}
\end{figure}

Figure~\ref{fig:ldm-vs-pinn-l01} and Figure~\ref{fig:ldm-vs-pinn-l04} show representative inverse modeling results for the two extremes of spatial correlation length: $\lambda = 0.1$ (high heterogeneity) and $\lambda = 0.4$ (smooth fields). 
In both cases, LD-DIM significantly outperforms the PINN-based method. 
At $\lambda = 0.1$, LD-DIM achieves a relative conductivity error of $\epsilon_K = 4.87 \times 10^{-1}$ and an SSIM of 0.709, while PINN yields a much larger error $\epsilon_K = 1.064$ and fails to reconstruct meaningful structure (SSIM = -0.054). 
For the smooth case $\lambda = 0.4$, the gap remains substantial: LD-DIM achieves $\epsilon_K = 3.193 \times 10^{-2}$ and SSIM = 0.987, compared to $\epsilon_K = 9.793 \times 10^{-1}$ and SSIM = 0.208 for PINN.

These results demonstrate that LD-DIM consistently achieves higher inversion accuracy by effectively capturing essential spatial heterogeneity, while avoiding the over-smoothed reconstructions commonly produced by PINN. The PINN approach, lacking explicit supervision on the conductivity field, struggles with high-dimensional optimization and does not possess inductive bias to reliably constrain the solution.

The superior performance of LD-DIM stems from its core design: it constrains the inverse problem within a low-dimensional, geologically informed latent space learned by the pretrained diffusion model. 
This constraint inherently steers the generative process toward physically plausible conductivity patterns, even under sparse supervision. 
Additionally, the differentiable forward model enables efficient gradient-based updates via reverse-mode automatic differentiation. 
In contrast, the PINN method relies solely on soft enforcement of PDE residuals and generic network expressiveness, which prove inadequate in severely ill-posed inverse problem.

\subsubsection{Incorporating conductivity measurements}\label{sec:including_k_data}
% --- Figure: λ = 0.1 ---
\begin{figure}[htbp!]
  \centering
  \includegraphics[width=0.8\textwidth]{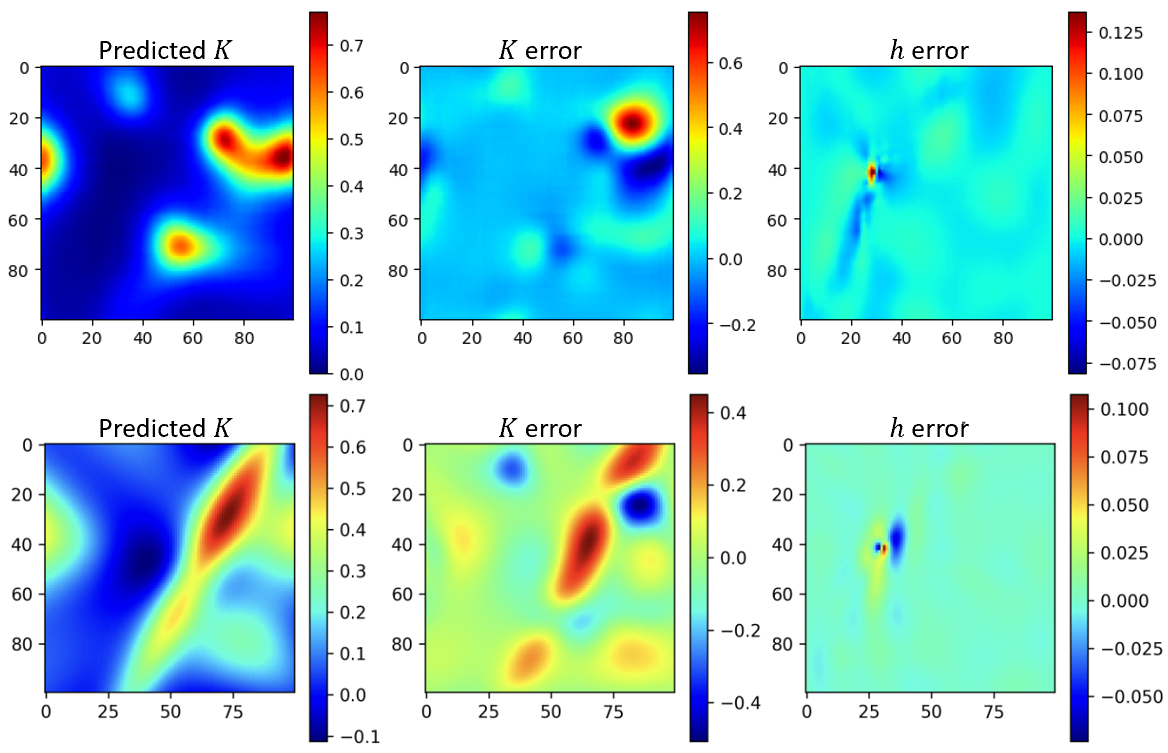}
  \caption{Comparison of LD-DIM and PINN at correlation length $\lambda=0.1$ with $5\times5$ conductivity observations.
  Top row: LD-DIM ($\epsilon_{K}=4.597\times10^{-1}$, SSIM = 0.767);
  bottom row: PINN ($\epsilon_{K}=5.388\times10^{-1}$, SSIM = 0.466).
  Columns (left to right): predicted conductivity field, error in conductivity, and error in hydraulic head.
  LD-DIM achieves superior accuracy and structural fidelity.}
  \label{fig:ldm-pinn-l01-K}
\end{figure}

% --- Figure: λ = 0.4 ---
\begin{figure}[htbp!]
  \centering
  \includegraphics[width=0.8\textwidth]{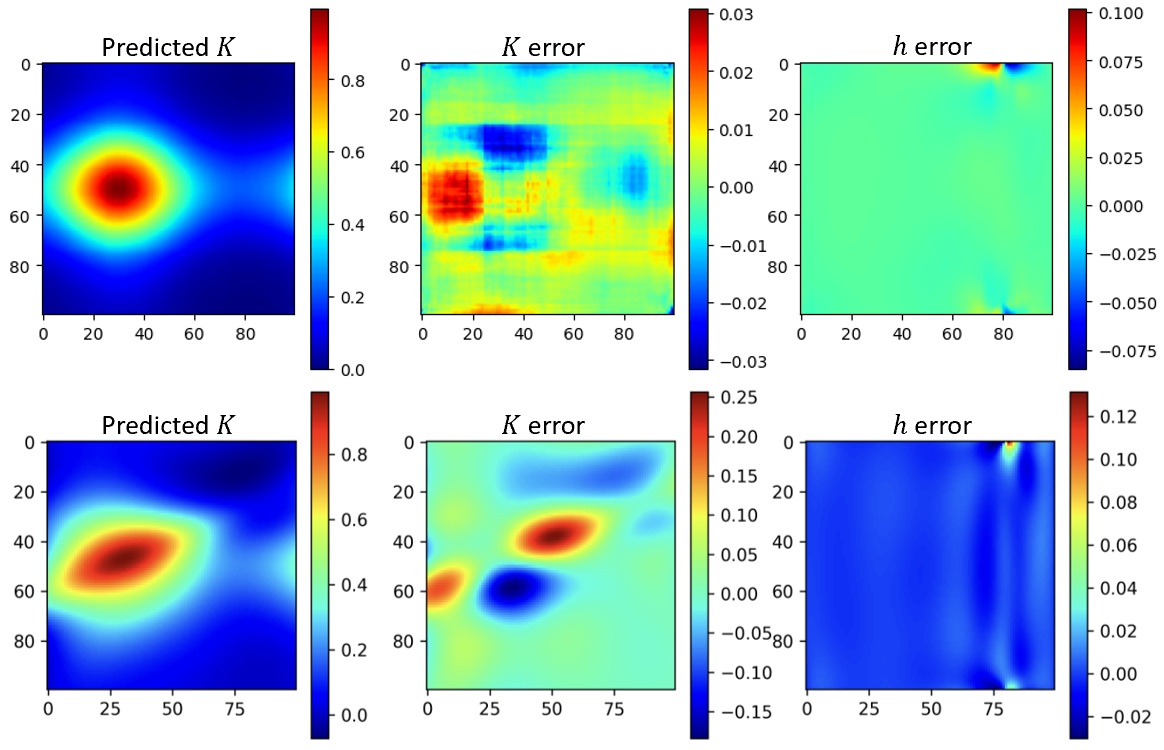}
  \caption{Comparison of LD-DIM and PINN at correlation length $\lambda=0.4$ with $5\times5$ conductivity observations.
  Top row: LD-DIM ($\epsilon_{K}=2.634\times10^{-2}$, SSIM = 0.993);
  bottom row: PINN ($\epsilon_{K}=1.882\times10^{-1}$, SSIM = 0.799).
  Columns (left to right): predicted conductivity field, error in conductivity, and error in hydraulic head.
  LD-DIM achieves superior accuracy and structural fidelity.}
  \label{fig:ldm-pinn-l04-K}
\end{figure}

In addition to the 256 hydraulic head observations, 25 conductivity measurements arranged in a $5\times5$ grid are incorporated into the inverse formulation.
Figs.~\ref{fig:ldm-pinn-l01-K} and~\ref{fig:ldm-pinn-l04-K} show that including this sparse conductivity information significantly improves the inversion accuracy for both approaches, especially for PINN. 
Compared with  cases without conductivity supervision (Figs.~\ref{fig:ldm-vs-pinn-l01} and~\ref{fig:ldm-vs-pinn-l04}), the relative error in PINN decreases from $1.064$ to $0.539$ for $\lambda=0.1$ and from $0.979$ to $0.188$ for $\lambda=0.4$, 
while the corresponding SSIM improves from $-0.054$ to $0.466$ and $0.208$ to $0.799$, respectively. 
These results indicate that even sparse conductivity data serve as informative priors that effectively enhance PINN optimization.

Consistent with the previous observations, when conductivity measurements are included, LD-DIM continues to outperform PINN, 
as evidenced by its higher estimation accuracy and structural fidelity:
$\epsilon_K = 0.460$ (SSIM = 0.767) and $\epsilon_K = 0.026$ (SSIM = 0.993) for $\lambda = 0.1$ and $\lambda = 0.4$, respectively.

These controlled comparisons in Sections \ref{subsec:pinn-comparison} and \ref{sec:including_k_data} demonstrate that LD-DIM’s advantage arises not merely from access to additional data, but from its ability to integrate prior knowledge through a learned latent space, which regularizes the inverse problem more effectively.
By encoding geological plausibility through the learned latent manifold, LD-DIM more effectively regularizes the inversion while preserving spatial heterogeneity, making it particularly valuable in practical settings where observations are limited.

Moreover, LD-DIM enforces physical constraints through a differentiable solver that satisfies the governing equations and boundary conditions exactly. 
This hard physical enforcement offers a substantially stronger constraint than the PINN-based methods, where physics is only weakly imposed through loss terms with empirically tuned weights.
As a result, LD-DIM avoids the need for manual hyperparameter tuning and delivers more robust performance across diverse spatial configurations.

%%%%%%%%%%%%%%%%%%%%%%%%%%%%%%%%%%%%
\subsection{Case study 2: Bimaterial conductivity fields}\label{subsec:bimaterial-fields}

This section extends the application of LD-DIM to bimaterial conductivity fields characterized by sharp material interfaces. 
Such patterns commonly arise in fractured aquifers, karst systems, and layered geological units, where strong hydraulic contrasts invalidate the smooth-conductivity assumptions typically used in hydrological modeling. 
Unlike the smoothly varying Gaussian fields in Section~\ref{subsec:gaussian-fields}, these discontinuous two-phase structures pose a significantly more challenging reconstruction problem, requiring accurate recovery of abrupt transitions and preservation of interface topology.

In this example, we first describe the generation of the bimaterial dataset (Section~\ref{subsec:bimaterial-generation}) and then assess the generative capacity of the latent diffusion model (Section~\ref{subsec:generative-validation}). 
Subsequent sections analyze inversion results under varying initialization seeds and  observation sparsity (Section~\ref{subsec:bimaterial-inversion}) and compare LD-DIM  against baseline methods (Section~\ref{subsec:vae-method-comparison}). 
Overall, these experiments evaluate LD-DIM’s ability to maintain geometric fidelity, recover interface morphology, and remain stable in heterogeneous two-phase systems.

\subsubsection{Synthetic data generation}\label{subsec:bimaterial-generation}

In this example, the bimaterial conductivity fields representing geologically complex scenarios with sharp permeability contrasts are constructed following the indicator-based stochastic modeling practices established in geostatistical reservoir characterization \cite{deutsch1992geostatistical}.  Each sample is formed by first constructing two log-normal fields representing the high- and low-conductivity materials,
\[
Y_1(\mathbf{x}) = \ln(10^{-5}) + \Phi_1(\mathbf{x}), \qquad
Y_2(\mathbf{x}) = \ln(10^{-8}) + \Phi_2(\mathbf{x}),
\]
where $\Phi_i$ are Gaussian processes with a Matérn covariance kernel ($\nu = 1$, $\lambda_1 = \lambda_2 = 50$). 

Then, a correlated splitting field $S(\mathbf{x})$ with correlation length $\lambda_{\mathrm{split}} = 200$ is used to assign the two components across the domain. The partition is determined by a quantile threshold,
\[
\alpha = 0.25 + 0.5\,p, \qquad p \sim \mathcal{U}(0,1),
\]
so that each material occupies between $25\%$ and $75\%$ of the domain. The spatial correlation in $S(\mathbf{x})$ leads to coherent two-phase patterns with sharp interfaces between two distinct conductivity regions.

%%%%%%%%%%%%%%%%%%%%%%%%%%%%%%%%%%%%%%%%%%%%%
\subsubsection{Generative model validation}\label{subsec:generative-validation}

Following the bimaterial data construction described above, a total of 2,000 conductivity-field samples are generated to capture a wide range of two-phase configurations and interface geometries. Among the 2,000 samples, 1,400 are used for training the generative prior, 400 for validation, and 200 for inverse modeling tests.

\begin{figure}[htbp!]
    \centering
    \begin{subfigure}[b]{\textwidth}
        \centering
        \includegraphics[width=0.65\textwidth]{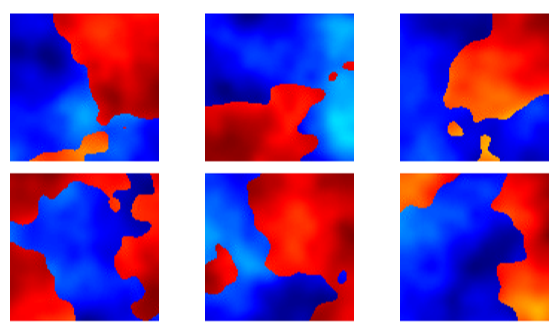}
        \caption{Ground truth bimaterial fields}
        \label{fig:ground-truth}
    \end{subfigure}

    \begin{subfigure}[b]{\textwidth}
        \centering
        \includegraphics[width=0.65\textwidth]{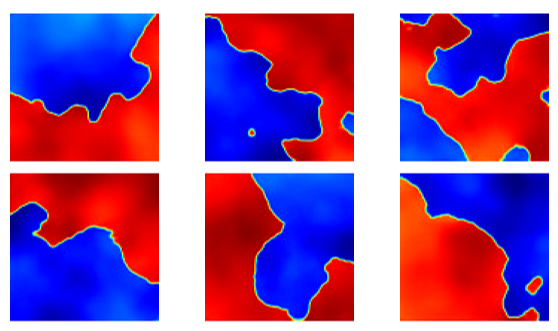}
        \caption{Latent diffusion model generated fields}
        \label{fig:ldm-generated}
    \end{subfigure}

    \begin{subfigure}[b]{\textwidth}
        \centering
        \includegraphics[width=0.65\textwidth]{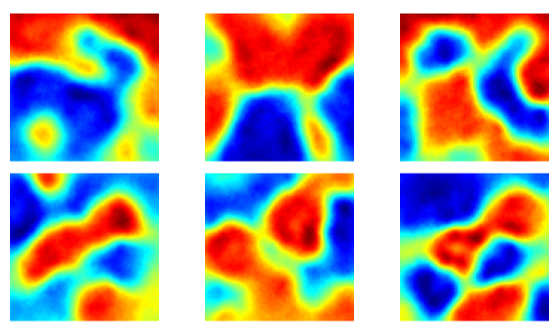}
        \caption{Variational autoencoder generated fields}
        \label{fig:vae-generated}
    \end{subfigure}
    
    \caption{Comparison of synthetic bimaterial conductivity fields generated by the latent diffusion model (LDM) and the baseline variational autoencoder (VAE). All fields are displayed in the $ln K$ domain.}
    \label{fig:generative-comparison}
\end{figure}

To assess the capability of the LDM generative prior, we compare it against a baseline variational autoencoder (VAE).
Figure~\ref{fig:generative-comparison} presents representative ground-truth bimaterial conductivity fields alongside conductivity samples generated by the LDM and the VAE. The results show that the LDM faithfully preserves sharp interfaces and the global two-phase structure, whereas the VAE produces overly smoothed transitions and blurred boundaries.

To quantify these differences, we apply two standard distributional metrics:
the Fréchet Inception Distance (FID)~\cite{heusel2017gans} and the Kernel Inception Distance (KID)~\cite{binkowski2018demystifying}, defined in Eqs. \eqref{eq:FID} and \eqref{eq:KID}, respectively. Both metrics operate on feature embeddings extracted from the edge representations of the fields. For each field $I$, the Sobel–Feldman edge detector $\mathcal{S}(\cdot)$ is used to compute the gradient magnitudes~\cite{sobel19683x3}, which indicate steep material jumps and suppresses smooth interior regions. 
These edge maps are normalized and encoded through a fixed feature extractor $\phi(\cdot)$, implemented as the pretrained Inception-based feature mapping used in standard FID/KID evaluations. The extractor is not trained within our pipeline; it provides a consistent embedding $\psi(I)=\phi(\mathcal{S}(I))$ that reflects the geometry, contrast, and spatial arrangement of interfaces.

\textbf{Fréchet Inception Distance (FID).}
Let $(\boldsymbol{\mu},\boldsymbol{\Sigma})$ and $(\boldsymbol{\mu}',\boldsymbol{\Sigma}')$ denote the mean and covariance of the embeddings $\{\psi(x_i)\}$ and $\{\psi(y_j)\}$ extracted from the real and generated fields. The distance is
\begin{equation}
\label{eq:FID}
\mathrm{FID}=
\bigl\lVert \boldsymbol{\mu}-\boldsymbol{\mu}' \bigr\rVert_2^2
+
\mathrm{tr}\!\left(
\boldsymbol{\Sigma}+\boldsymbol{\Sigma}'
-2\,(\boldsymbol{\Sigma}\boldsymbol{\Sigma}')^{1/2}
\right),
\end{equation}
which measures discrepancies in the global statistics of the embedded edge features. Lower FID values indicate closer agreement in interface geometry and contrast.

\textbf{Kernel Inception Distance (KID).}
Given feature sets $\{\psi(x_i)\}_{i=1}^{N}$ and $\{\psi(y_j)\}_{j=1}^{M}$ from real and generated fields, the unbiased estimator is
\begin{equation}
\label{eq:KID}
\begin{aligned}
\mathrm{KID} &=
\frac{1}{N(N-1)}
\sum_{\substack{i,j=1\\ i\neq j}}^{N}
k\!\big(\psi(x_i),\psi(x_j)\big)
+
\frac{1}{M(M-1)}
\sum_{\substack{i,j=1\\ i\neq j}}^{M}
k\!\big(\psi(y_i),\psi(y_j)\big)
\\[3pt]
&\quad-
\frac{2}{NM}
\sum_{i=1}^{N}\sum_{j=1}^{M}
k\!\big(\psi(x_i),\psi(y_j)\big).
\end{aligned}
\end{equation}
where $k(u,v)=\left(\frac{1}{d}u^\top v+1\right)^3$. Because KID directly compares pairwise similarities in the embedded edge features, lower values indicate a better match to the reference distribution in terms of interface sharpness and fine-scale structural patterns.

The numerical values of these metrics, summarized in Table~\ref{tab:generative-metrics}, further reinforce the qualitative observations above: the LDM achieves an FID more than five times lower than the VAE, and its KID is reduced by over fourfold. These substantially lower scores indicate a closer match to the reference distribution in both local boundary sharpness and overall interface organization, demonstrating that the LDM is well suited for generating physically realistic bimaterial fields to serve as priors in inverse modeling.

\begin{table}[htbp!]
\centering
\caption{Quantitative evaluation of generative performance for the latent diffusion model (LDM) and the variational autoencoder (VAE) using FID and KID. Lower scores indicate better preservation of sharp interfaces and structural details.}
\vspace{0.3cm}
\begin{tabular}{@{}lcc@{}}
\toprule
\textbf{Metric} & \textbf{Latent Diffusion Model} & \textbf{Variational Autoencoder} \\
\midrule
FID ($\downarrow$) & 77.4 & 487.1 \\
KID ($\downarrow$) & 0.1102 & 0.5452 \\
\bottomrule
\end{tabular}
\label{tab:generative-metrics}
\end{table}

%%%%%%%%%%%%%%%%%%%%%%%%%%%%%%%%%%%%%%
\subsubsection{Inversion performance analysis}\label{subsec:bimaterial-inversion}
We evaluate the inversion capability of LD-DIM on bimaterial conductivity fields characterized by sharp interfaces and strong conductivity contrasts.
 Two key aspects are investigated:
%Specifically, we examine two key aspects: 
(i) sensitivity to random initialization and (ii) the effect of observation density on reconstruction quality and stability.
In addition, we compare LD-DIM with a physics-embedded VAE-based framework~\citep{wu2023physics}, demonstrating that LD-DIM produces sharper interface representations and consistently higher reconstruction fidelity.

As in the Gaussian case study presented in Section~\ref{subsec:gaussian-fields}, all bimaterial inversion experiments are conducted on sample fields excluded from the training dataset, ensuring an unbiased assessment of generalization performance.

%%%%%%%%%%%%%%%%%%%%%
\subsubsection*{(i) Sensitivity analysis of initial seed}\label{para:seed-sensitivity}
The inherent ill-posedness of the inverse problem introduces variability in reconstruction outcomes, making it important to examine how different random initialization seeds affect the predicted fields. To assess this effect, we consider three representative realizations obtained from different seed initializations for a bimaterial test case using 25 head observations on a $5\times5$ uniform grid.

\begin{figure}[H]
    \centering
    \includegraphics[width=1.0\textwidth]{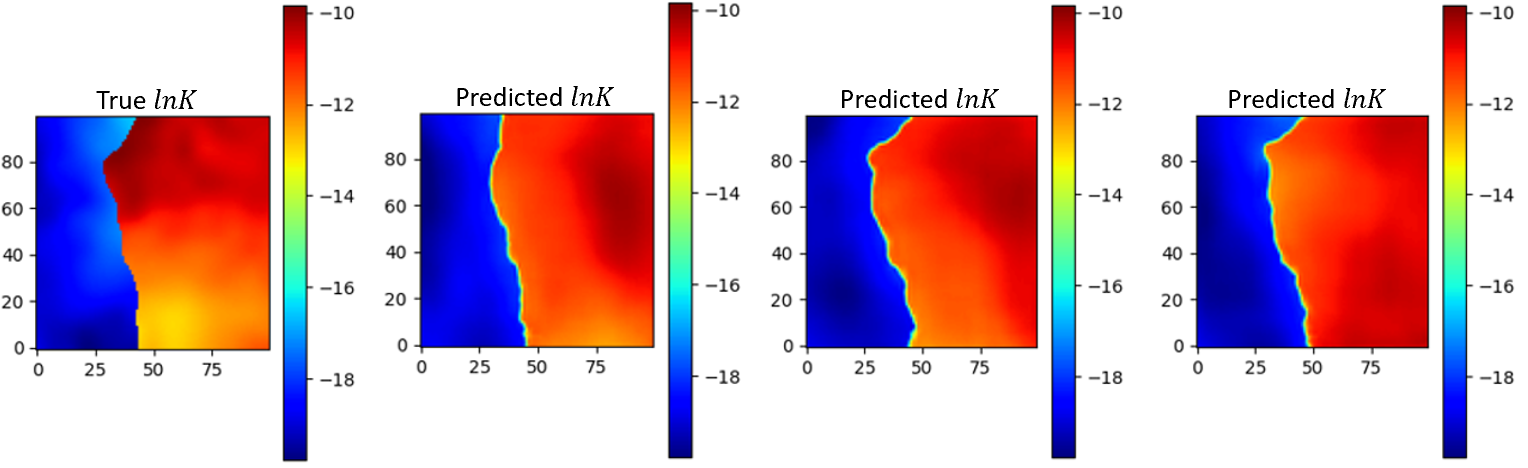}
    \caption{Representative examples of sensitivity analysis for bimaterial conductivity field reconstruction using 25 observations with 3 random initialization seeds. From left to right: true log-conductivity field, and three representative predictions with relative errors of 0.275 (SSIM = 0.818), 0.301 (SSIM = 0.807), and 0.329 (SSIM = 0.773). (All fields are displayed and evaluated in the $ln K$ domain.)}
    \label{fig:bimaterial-sparse}
\end{figure}

Fig.~\ref{fig:bimaterial-sparse} shows conductivity reconstructions from three different random initialization seeds.   Despite starting from different initial latent vectors, all three realizations successfully recover the sharp two-phase structure and major geometric features. The mean-corrected relative errors $\tilde{\epsilon}_{K}$ associated with the 3 different initializations remain within a narrow range (0.275, 0.301, and 0.329), with corresponding SSIM values of 0.818, 0.807, and 0.773, respectively. Only minor variations in interface location are observed across seeds.  
This robustness across initializations can be attributed to two key mechanisms: first, the latent diffusion prior constrains the solution space to geologically plausible bimaterial realizations, preventing convergence to physically unrealistic field patterns; second, the explicit regularization term in Eq. \eqref{eq:chain_rule} maintains proximity to the learned latent distribution, further stabilizing the optimization process.

\begin{figure}[H]
    \centering
    \includegraphics[width=\textwidth]{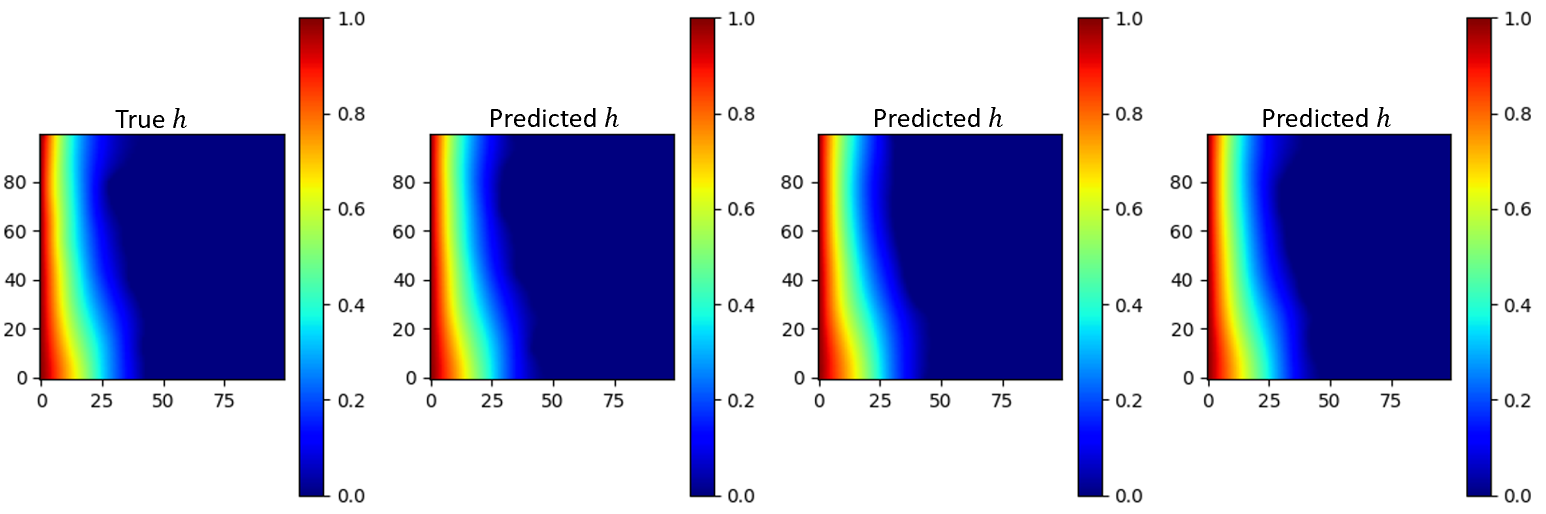}
    \caption{Sensitivity of hydraulic head field reconstruction (25 observations) to 3 random initialization seeds. From left to right: ground truth hydraulic head field and three representative predictions from different random seeds. The corresponding metrics are: relative error $=2.856\times10^{-2}$, SSIM $=0.992$; relative error $=3.241\times10^{-2}$, SSIM $=0.982$; relative error $=3.074\times10^{-2}$, SSIM $=0.987$.}
    \label{fig:head-seed-sensitivity}
\end{figure}

The hydraulic head reconstructions are shown in Fig.~\ref{fig:head-seed-sensitivity}, exhibiting even lower sensitivity to initial seeds.  The three realizations produce relative errors of $2.856\times10^{-2}$, $3.241\times10^{-2}$, and $3.074\times10^{-2}$, with SSIM values of 0.992, 0.982, and 0.987. 
The predicted head fields differ only slightly across seeds and remain very close to the reference solution. 
This consistency indicates that hydraulic head fields are substantially easier to reconstruct than the high-contrast conductivity field, due to both the direct availability of head measurements and the inherently smoother nature of the governing PDE solution. As a result, even when conductivity reconstructions exhibit geometric variations near material interfaces, the corresponding head predictions show very limited seed-induced variability.

\subsubsection*{(ii) Observation-density analysis}\label{para:Observation-density}

\begin{figure}[htbp!]
    \centering
    \includegraphics[width=\textwidth]{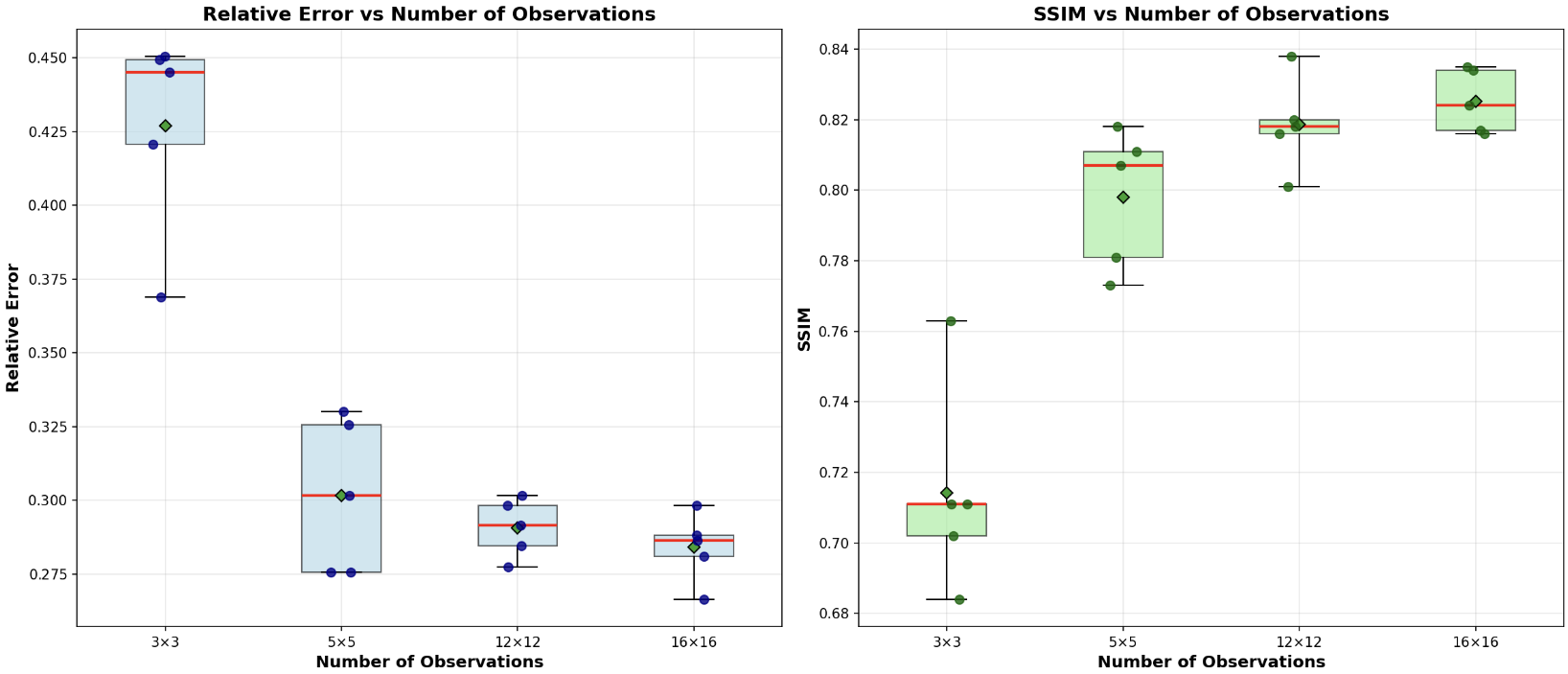}
    \caption{Analysis of observation density on log-conductivity field reconstruction quality. 
    Left: relative error vs.\ number of observations; Right: SSIM vs.\ number of observations. 
    Each box plot summarizes results across multiple random seeds for each layout (3$\times$3, 5$\times$5, 12$\times$12, 16$\times$16). 
    The red horizontal line indicates the median, and the diamond marker denotes the mean value computed over all random seeds. (All fields are displayed and evaluated in the $ln K$ domain.)}
    \label{fig:obs-density-K}
\end{figure}

\begin{figure}[H]
    \centering
    \includegraphics[width=\textwidth]{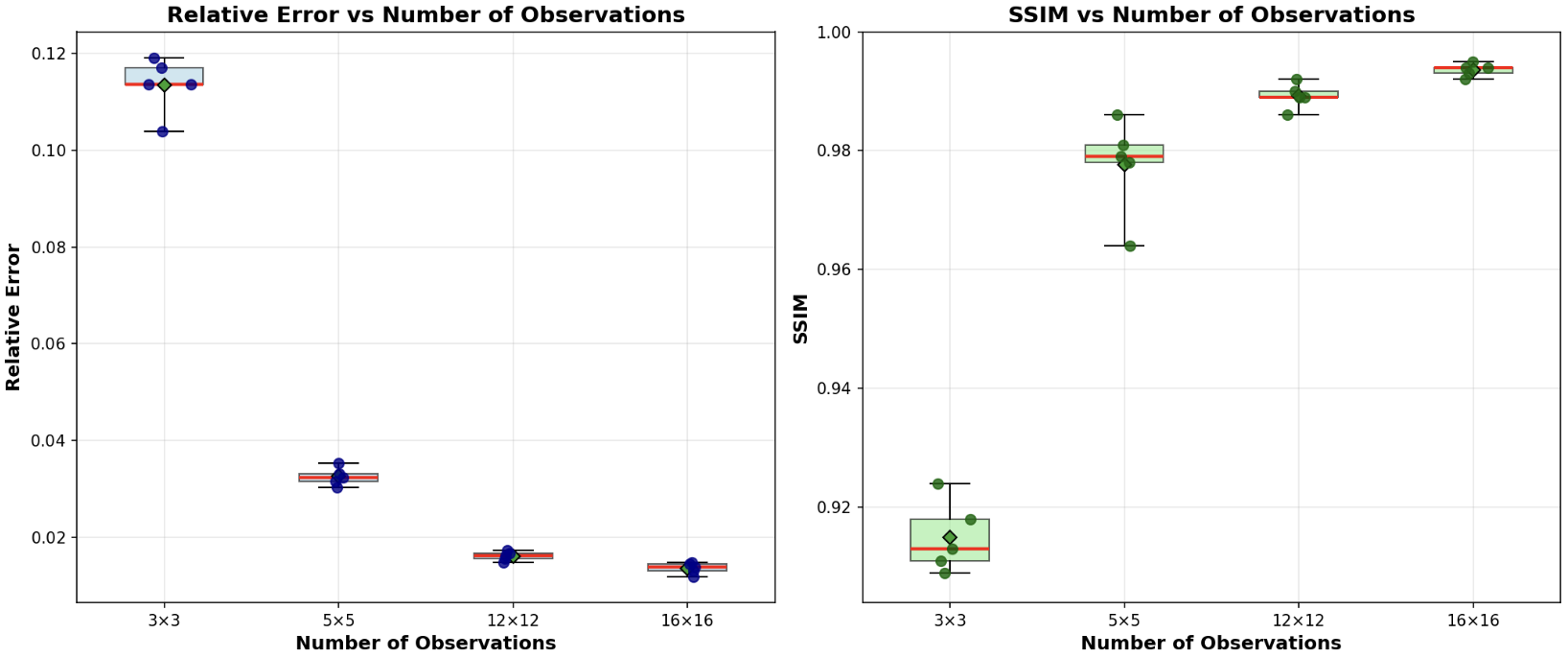}
    \caption{Analysis of observation density on hydraulic head field reconstruction quality. 
    Left: relative error vs.\ number of observations; Right: SSIM vs.\ number of observations. 
    Each box plot summarizes results across multiple random seeds for each layout (3$\times$3, 5$\times$5, 12$\times$12, 16$\times$16). 
    The red horizontal line indicates the median, and the diamond marker denotes the mean value computed over all random seeds.}
    \label{fig:obs-density-h}
\end{figure}

In practical field applications, the number of available head measurements is often limited. To understand how observation density influences reconstruction performance, we systematically evaluate the framework using four uniform observation grids: $3\times3$, $5\times5$, $12\times12$, and $16\times16$ grids. This analysis quantifies the trade-off between data availability and reconstruction fidelity for both conductivity and hydraulic head fields, as shown in Fig.~\ref{fig:obs-density-K} and Fig.~\ref{fig:obs-density-h}.

For the conductivity reconstruction (Fig.~\ref{fig:obs-density-K}), the $3\times3$ grid results in a median relative error of approximately $0.44$ and a median SSIM of $0.71$, indicating that $9$ head observations provide only limited constraint on the conductivity field reconstruction. Increasing the density to a $5\times5$ grid reduces the median error to about $0.30$ and raises the SSIM to roughly $0.81$, reflecting improved recovery of the conductivity field. 
The $12\times12$ and $16\times16$ grids further enhance the consistency across different seeds, with the latter achieving a median error of $0.28$ and a median SSIM of $0.82$, representing the highest interface fidelity of the tested configurations.
The progressive narrowing of the interquartile ranges across the four densities demonstrates reduced sensitivity to initialization and increased robustness of the reconstructed conductivity fields. 
This robustness stems from the inherent regularization of the LD-DIM framework: by constraining  inversion within a geologically informed latent manifold, the diffusion prior compensates for sparse observations and favors structurally plausible conductivity patterns. 
Consequently, LD-DIM delivers stable and reliable reconstructions even under limited data availability.

The hydraulic head reconstruction results (Fig.~\ref{fig:obs-density-h}) further illustrate the influence of the number of head observations on both the accuracy and the stability of the reconstructed head fields. 
The median relative error decreases sharply from $0.11$ to $0.03$ as the observation grid increases from a $3 \times 3$ to a $5 \times 5$, showing that even a modest increase in observation density substantially constrains the groundwater flow solution. 
Beyond $5 \times 5$, additional observations provide diminishing returns in absolute accuracy, e.g., the SSIM is already close to $0.99$, but they do meaningfully reduce the variance across different initial seeds.

%%%%%%%%%%%%%%%%%%%%%%%%%%%%%%%%%%%%%%%
\subsubsection*{(iii) Comparison with physics-embedded VAE method}\label{subsec:vae-method-comparison}

We compare LD-DIM with the physics-embedded VAE approach of \cite{wu2023physics}, which employs a variational autoencoder trained on bimaterial fields as a generative prior. During inversion, the latent variables are optimized to minimize discrepancy between simulated and observed hydraulic heads while maintaining proximity to the learned latent distribution.

\begin{figure}[H]
    \centering
    \includegraphics[width=1.0\textwidth]{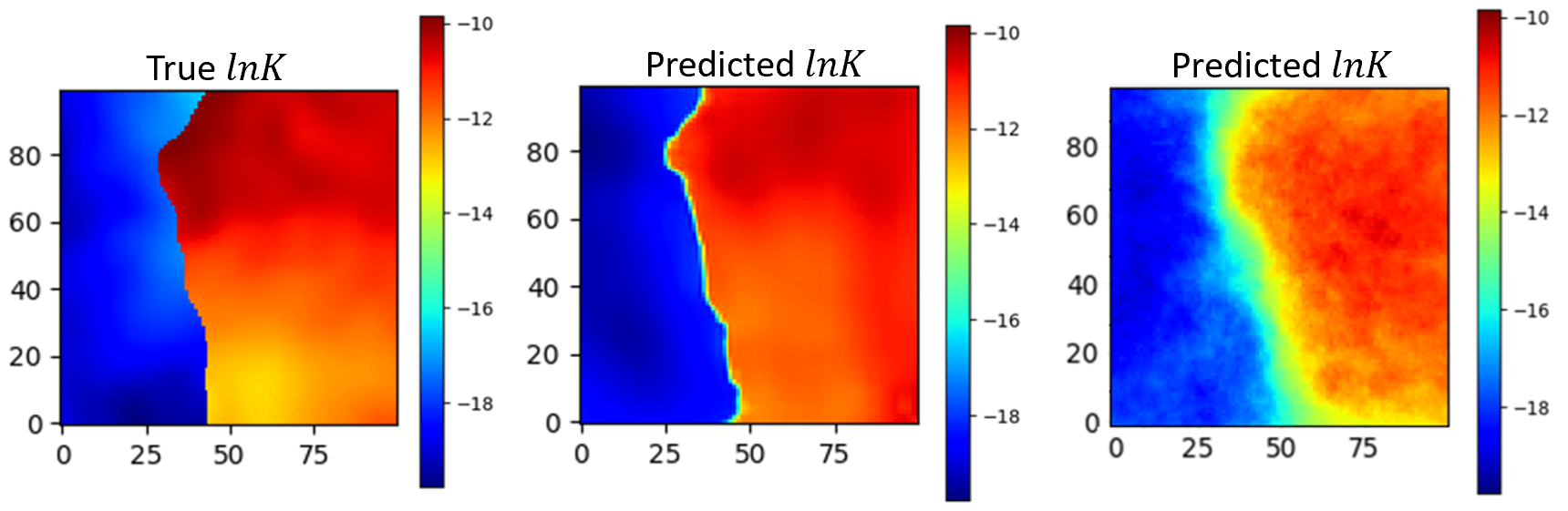}
    \caption{Comparative reconstruction of a representative bimaterial field. Left: true field; center: LD-DIM result ($\tilde{\epsilon}_{K}$ = 0.272); right: physics-embedded VAE result ($\tilde{\epsilon}_{K} = 0.387$). LD-DIM achieves superior interface preservation and reduced reconstruction error. (All fields are displayed and evaluated in the $ln K$ domain.)}
    \label{fig:vae-comparison}
\end{figure}

The comparison in Figure~\ref{fig:vae-comparison} highlights the advantages of LD-DIM for bimaterial reconstruction. 
LD-DIM reduces the relative error $\tilde{\epsilon}_{K}$ from $0.387$ to $0.272$ and achieves a higher SSIM, reflecting improved recovery of both global structure and local boundary detail. While The VAE-based method is able to capture the overall bimaterial interface geometry,  
it produces noticeably smoothed boundaries and fails to preserve the sharp transitions characteristic of high-contrast fields.
This limitation stems from the Gaussian prior commonly used in VAEs, which biases reconstructions toward smooth, low-frequency patterns and struggles to represent discontinuities.

In contrast, the latent diffusion model leverages an iterative denoising process that approximates a richer, non-Gaussian latent distribution and progressively refines high-frequency spatial features essential for accurate sharp-interface representations. 
As a result, LD-DIM provides a more expressive and geologically realistic prior, leading to substantially improved interface fidelity and more accurate reconstruction of flow-controlling structures.
Overall, the successful reconstruction of complex sharp-interface structures underscores the potential of LD-DIM for practical applications in fractured aquifers, layered subsurface systems, and other geologically heterogeneous environments.

\section{Conclusion}\label{sec:conclusion}

This study introduced LD-DIM, a latent diffusion–based differentiable inversion framework that integrates a physics-informed generative prior with a differentiable finite-volume solver for subsurface inverse modeling. 
By embedding inversion within a low-dimensional latent space while enforcing physical consistency through exact numerical discretization, LD-DIM effectively mitigates the ill-posedness inherent in high-dimensional subsurface inverse problems. 
The effectiveness of the proposed LD-DIM framework arises from three key features. First, latent-space regularization confines optimization to a compact yet geologically plausible manifold, alleviating the curse of dimensionality that limits conventional gradient-based inversion in high-dimensional parameter spaces. 
Second, the latent diffusion prior preserves both smooth spatial correlations and sharp geological discontinuities, overcoming common limitations of existing approaches, including interface blurring in VAE-based methods, training instability in GAN-based models, and insufficient regularization in PINN-based inversions. 
%Third, adjoint-based exact gradient computation enables efficient end-to-end optimization with transparent sensitivity analysis, avoiding the substantial computational cost associated with ensemble-based calibration and Bayesian sampling methods while maintaining strict physical consistency.
Third, the proposed hybrid adjoint–automatic differentiation strategy enables exact gradient computation for PDE-constrained inverse problems, supporting efficient end-to-end optimization with transparent sensitivity analysis while maintaining strict physical consistency.
%while avoiding the substantial computational cost associated with ensemble-based calibration and Bayesian sampling methods.
% The strong performance of LD-DIM in this challenging setting can be attributed to three key factors: (i) the geometric prior embedded in the latent diffusion model preserves sharp geological discontinuities, enabling accurate interface recovery; (ii) the low-dimensional latent representation yields a tractable optimization landscape even for complex geometries; and (iii) the differentiable physics solver enforces exact satisfaction of the governing groundwater flow equations, ensuring physically coherent predictions.

We evaluate the proposed LD-DIM framework using two representative case studies involving Gaussian random fields and bimaterial conductivity fields, and systematically examine its generative capability as well as its sensitivity to measurement availability and initialization.
For Gaussian conductivity fields with varying correlation lengths ($\lambda = 0.1, 0.2, 0.3, 0.4$), 
LD-DIM demonstrates consistently strong performance, achieving low reconstruction error and high structural similarity even in the most challenging low-correlation regime. 
Compared with the PINN-based inversion baseline, LD-DIM reduces reconstruction error by approximately $2.2\times$ at $\lambda = 0.1$ (from $\epsilon_K = 1.064$ to $0.487$) and by more than $30\times$ at $\lambda = 0.4$ (from $\epsilon_K = 0.9793$ to $0.0319$), while consistently preserving fine-scale heterogeneity.

For bimaterial conductivity fields, LD-DIM effectively resolved sharp geological discontinuities. 
Relative to a baseline VAE, the latent diffusion prior exhibits superior generative expressiveness, and the resulting inversion achieves substantially lower reconstruction error ($\epsilon_K = 0.272$ vs.\ $0.387$) while accurately recovering material interfaces critical to fractured and layered subsurface systems.
Comprehensive sensitivity analyses across different observation densities and random initialization seeds further demonstrate the robustness of the proposed framework. 
Statistical analyses of reconstructed conductivity and hydraulic head fields reveal stable medians and decreasing interquartile ranges as observation density increases, confirming the consistency and reliability of LD-DIM.

Looking ahead, the proposed LD-DIM framework is readily extensible to a broad class of physics-constrained inverse problems beyond subsurface flow, including seismic full-waveform inversion for  velocity characterization, electromagnetic inversion for mineral exploration, thermal modeling in geothermal systems, and multiphase flow inversion in energy and environmental applications. These results underscore the broader potential of physics-constrained latent diffusion models as a scalable and computationally efficient paradigm for inverse modeling in complex geophysical and engineering systems.

\section*{Data availability}\label{Data} 
The datasets used in this study will be made publicly available at \url{https://github.com/IntelligentMechanicsLab/LD-DIM}.

% --------------------------------------------------------------------
%%%%%%%%%%%
\appendix
\section{Finite volume discretization}\label{sec:app_A}

This appendix provides additional details on the finite volume method (FVM) used in Section~\ref{sec:fvm} to discretize the steady-state groundwater flow equation.

We consider a 2D domain $\Omega$ discretized using a uniform Cartesian mesh with cell-centered variables and spacings $\Delta x = \Delta y$. Let $\hat{h}_{i,j}$ denote the pressure head at the center of grid cell $(i,j)$. Integrating the governing equation~\eqref{eq:darcy_governing} over a control volume and applying the divergence theorem yields a discrete flux balance:
\begin{equation}
F_{i+\frac12,j} - F_{i-\frac12,j} + F_{i,j+\frac12} - F_{i,j-\frac12} = 0 ,
\label{eq:cv-balance}
\end{equation}
where $F_{i\pm\frac12,j}$ and $F_{i,j\pm\frac12}$ denote diffusive fluxes through the four faces of cell.

Each face flux is approximated via two-point flux approximation (TPFA) as the product of a transmissibility coefficient and the head difference across adjacent cells:
\begin{align}
F_{i+\frac12,j} &= T_{i+\frac12,j}\,\big(\hat{h}_{i+1,j}-\hat{h}_{i,j}\big), \quad
F_{i-\frac12,j} = T_{i-\frac12,j}\,\big(\hat{h}_{i,j}-\hat{h}_{i-1,j}\big), \\
F_{i,j+\frac12} &= T_{i,j+\frac12}\,\big(\hat{h}_{i,j+1}-\hat{h}_{i,j}\big), \quad
F_{i,j-\frac12} = T_{i,j-\frac12}\,\big(\hat{h}_{i,j}-\hat{h}_{i,j-1}\big),
\end{align}
which maintains consistency on $K$-orthogonal grids for isotropic media \citep{aavatsmark2002introduction,eymard2000finite}.
For heterogeneous isotropic media, the face-normal conductivity is computed as the harmonic mean of
adjacent cell values,
\begin{align}
\hat{K}^{\rm harm}_{i+\frac12,j} &= \frac{2\,\hat{K}_{i,j}\,\hat{K}_{i+1,j}}{\hat{K}_{i,j}+\hat{K}_{i+1,j}}, \quad
\hat{K}^{\rm harm}_{i-\frac12,j} = \frac{2\,\hat{K}_{i,j}\,\hat{K}_{i-1,j}}{\hat{K}_{i,j}+\hat{K}_{i-1,j}}, \\
\hat{K}^{\rm harm}_{i,j+\frac12} &= \frac{2\,\hat{K}_{i,j}\,\hat{K}_{i,j+1}}{\hat{K}_{i,j}+\hat{K}_{i,j+1}}, \quad
\hat{K}^{\rm harm}_{i,j-\frac12} = \frac{2\,\hat{K}_{i,j}\,\hat{K}_{i,j-1}}{\hat{K}_{i,j}+\hat{K}_{i,j-1}},
\end{align}
such that the per-face transmissibilities on a uniform grid become
\begin{equation}
T_{i\pm\frac12,j}=\hat{K}^{\rm harm}_{i\pm\frac12,j}\frac{\Delta y}{\Delta x},
\qquad
T_{i,j\pm\frac12}=\hat{K}^{\rm harm}_{i,j\pm\frac12}\frac{\Delta x}{\Delta y},
\end{equation}
following standard finite-volume diffusive conductance formulas \citep{patankar2018numerical,eymard2000finite}.
Substituting the face fluxes into Equation~\eqref{eq:cv-balance} produces a sparse five-point stencil. 
Using the notation $p\equiv(i,j)$ and $E=(i+1,j)$, $W=(i-1,j)$, $N=(i,j+1)$, $S=(i,j-1)$, the matrix entries are
\begin{align}
A_{p,p} &= -\!\left(
T_{i+\frac12,j}+T_{i-\frac12,j}+T_{i,j+\frac12}+T_{i,j-\frac12}
\right), \\
A_{p,E} &= T_{i+\frac12,j},\quad
A_{p,W} = T_{i-\frac12,j},\quad
A_{p,N} = T_{i,j+\frac12},\quad
A_{p,S} = T_{i,j-\frac12},
\end{align}
and assembling over all cells gives the linear system
\begin{equation}
% A(\hat{K})\,\hat{h} = b.
\bm{A}(\hat{\bm K})\,\hat{\bm{h}}=\bm{b},
\end{equation}
where $\bm{A}$ is the stiffness matrix dependent on the discrete conductivity field $\hat{\bm K}$, and $\bm{b}$ encodes the boundary conditions.
As described in Section \ref{sec:fvm},
the Dirichlet and Neumann boundary conditions can be imposed on this linear system accordingly.

%%%%%%%%%
\section{Numerical Implementation of the FVM Solver}
\label{app:implementation}

The sparse linear system in \eqref{eq:fvm_solv}, $\bm{A}\hat{\bm{h}} = \bm{b}$, which results from the finite-volume discretization described in the main text, is constructed and solved using computationally efficient techniques. The sparse matrix $\bm{A}$ is first assembled in the \emph{list-of-lists} (LIL) format, which is well-suited for the incremental insertion of stencil coefficients during the iterative assembly process over the grid cells. Once the matrix coefficients are fully populated, the matrix is converted to the \emph{compressed sparse row} (CSR) format, which is highly optimized for the fast matrix-vector operations required by numerical solvers. The linear system is then solved using the \texttt{spsolve} function from the SciPy library, a high-performance direct sparse solver based on the SuperLU library.

To enable gradient-based optimization for the inverse problem, the entire forward solver is implemented to be end-to-end differentiable. This is achieved by defining a custom \emph{vector–Jacobian product} (VJP) in the JAX framework, which provides the core mechanism for reverse-mode automatic differentiation. For our solver, the custom VJP takes the gradient of the loss function $\ell$ with respect to the hydraulic head $\hat{\bm{h}}$ (i.e., $\nabla_{\hat{\bm{h}}}\ell$), and constructs the corresponding adjoint system. 

The adjoint system is obtained by modifying the forward operator to handle the appropriate boundary conditions: at Dirichlet boundaries where the forward problem prescribes head values, the adjoint field $\bm{\lambda}$ is constrained to zero. The modified adjoint system
\begin{equation}
\bm{A}_{\text{adj}}\bm{\lambda} = \nabla_{\hat{\bm{h}}}\ell
\end{equation}
is then solved to obtain the adjoint field $\bm{\lambda}$, where $\bm{A}_{\text{adj}}$ incorporates the proper adjoint boundary conditions. This enables the efficient computation of gradients with respect to the conductivity field, which are subsequently backpropagated through the generative model via the chain rule. This approach ensures that the entire finite-volume operator is seamlessly embedded within the end-to-end optimization framework.

\bibliographystyle{apalike}

\bibliography{ref_subsurface,ref_diffusion}

@article{dasgupta2025conditional,
  title={Conditional score-based diffusion models for solving inverse elasticity problems},
  author={Dasgupta, Agnimitra and Ramaswamy, Harisankar and Murgoitio-Esandi, Javier and Foo, Ken Y and Li, Runze and Zhou, Qifa and Kennedy, Brendan F and Oberai, Assad A},
  journal={Computer Methods in Applied Mechanics and Engineering},
  volume={433},
  pages={117425},
  year={2025},
  publisher={Elsevier}
}

@article{bezgin2023jax,
  title={JAX-Fluids: A fully-differentiable high-order computational fluid dynamics solver for compressible two-phase flows},
  author={Bezgin, Deniz A and Buhendwa, Aaron B and Adams, Nikolaus A},
  journal={Computer Physics Communications},
  volume={282},
  pages={108527},
  year={2023},
  publisher={Elsevier}
}

@article{xue2023jax,
  title={JAX-FEM: A differentiable GPU-accelerated 3D finite element solver for automatic inverse design and mechanistic data science},
  author={Xue, Tianju and Liao, Shuheng and Gan, Zhengtao and Park, Chanwook and Xie, Xiaoyu and Liu, Wing Kam and Cao, Jian},
  journal={Computer Physics Communications},
  volume={291},
  pages={108802},
  year={2023},
  publisher={Elsevier}
}

@article{karniadakis2021physics,
  title={Physics-informed machine learning},
  author={Karniadakis, George Em and Kevrekidis, Ioannis G and Lu, Lu and Perdikaris, Paris and Wang, Sifan and Yang, Liu},
  journal={Nature Reviews Physics},
  volume={3},
  number={6},
  pages={422--440},
  year={2021},
  publisher={Nature Publishing Group UK London}
}

@article{xiao2021deep,
  title={Deep-learning-based adjoint state method: Methodology and preliminary application to inverse modeling},
  author={Xiao, Cong and Deng, Ya and Wang, Guangdong},
  journal={Water Resources Research},
  volume={57},
  number={2},
  pages={e2020WR027400},
  year={2021},
  publisher={Wiley Online Library}
}

@article{wang2004image,
  title={Image quality assessment: from error visibility to structural similarity},
  author={Wang, Zhou and Bovik, Alan C and Sheikh, Hamid R and Simoncelli, Eero P},
  journal={IEEE transactions on image processing},
  volume={13},
  number={4},
  pages={600--612},
  year={2004},
  publisher={IEEE}
}

@article{brunet2011mathematical,
  title={On the mathematical properties of the structural similarity index},
  author={Brunet, Dominique and Vrscay, Edward R and Wang, Zhou},
  journal={IEEE Transactions on Image Processing},
  volume={21},
  number={4},
  pages={1488--1499},
  year={2011},
  publisher={IEEE}
}

@article{he2020physics,
  title={Physics-informed neural networks for multiphysics data assimilation with application to subsurface transport},
  author={He, QiZhi and Barajas-Solano, David and Tartakovsky, Guzel and Tartakovsky, Alexandre M},
  journal={Advances in Water Resources},
  volume={141},
  pages={103610},
  year={2020},
  publisher={Elsevier}
}

@article{binkowski2018demystifying,
  title={Demystifying mmd gans},
  author={Bi{\'n}kowski, Miko{\l}aj and Sutherland, Danica J and Arbel, Michael and Gretton, Arthur},
  journal={arXiv preprint arXiv:1801.01401},
  year={2018}
}

@article{heusel2017gans,
  title={Gans trained by a two time-scale update rule converge to a local nash equilibrium},
  author={Heusel, Martin and Ramsauer, Hubert and Unterthiner, Thomas and Nessler, Bernhard and Hochreiter, Sepp},
  journal={Advances in neural information processing systems},
  volume={30},
  year={2017}
}

@article{neuman1973calibration,
  title={Calibration of distributed parameter groundwater flow models viewed as a multiple-objective decision process under uncertainty},
  author={Neuman, Shlomo P},
  journal={Water Resources Research},
  volume={9},
  number={4},
  pages={1006--1021},
  year={1973},
  publisher={Wiley Online Library}
}

@article{yeh1986review,
  title={Review of parameter identification procedures in groundwater hydrology: The inverse problem},
  author={Yeh, William W-G},
  journal={Water resources research},
  volume={22},
  number={2},
  pages={95--108},
  year={1986},
  publisher={Wiley Online Library}
}

@article{mclaughlin1996reassessment,
  title={A reassessment of the groundwater inverse problem},
  author={McLaughlin, Dennis and Townley, Lloyd R},
  journal={Water Resources Research},
  volume={32},
  number={5},
  pages={1131--1161},
  year={1996},
  publisher={Wiley Online Library}
}

@article{jafarpour2008history,
  title={History matching with an ensemble Kalman filter and discrete cosine parameterization},
  author={Jafarpour, Behnam and McLaughlin, Dennis B},
  journal={Computational Geosciences},
  volume={12},
  number={2},
  pages={227--244},
  year={2008},
  publisher={Springer}
}

@article{laloy2017inversion,
  title={Inversion using a new low-dimensional representation of complex binary geological media based on a deep neural network},
  author={Laloy, Eric and H{\'e}rault, Romain and Lee, John and Jacques, Diederik and Linde, Niklas},
  journal={Advances in water resources},
  volume={110},
  pages={387--405},
  year={2017},
  publisher={Elsevier}
}

@article{laloy2018training,
  title={Training-image based geostatistical inversion using a spatial generative adversarial neural network},
  author={Laloy, Eric and H{\'e}rault, Romain and Jacques, Diederik and Linde, Niklas},
  journal={Water Resources Research},
  volume={54},
  number={1},
  pages={381--406},
  year={2018},
  publisher={Wiley Online Library}
}

@book{ramm2005inverse,
  title={Inverse problems: mathematical and analytical techniques with applications to engineering},
  author={Ramm, Alexander G},
  year={2005},
  publisher={Springer}
}

@article{linde2015geological,
  title={Geological realism in hydrogeological and geophysical inverse modeling: A review},
  author={Linde, Niklas and Renard, Philippe and Mukerji, Tapan and Caers, Jef},
  journal={Advances in Water Resources},
  volume={86},
  pages={86--101},
  year={2015},
  publisher={Elsevier}
}

@article{zhou2014inverse,
  title={Inverse methods in hydrogeology: Evolution and recent trends},
  author={Zhou, Haiyan and G{\'o}mez-Hern{\'a}ndez, J Jaime and Li, Liangping},
  journal={Advances in Water Resources},
  volume={63},
  pages={22--37},
  year={2014},
  publisher={Elsevier}
}

@article{sarma2008kernel,
  title={Kernel principal component analysis for efficient, differentiable parameterization of multipoint geostatistics},
  author={Sarma, Pallav and Durlofsky, Louis J and Aziz, Khalid},
  journal={Mathematical Geosciences},
  volume={40},
  number={1},
  pages={3--32},
  year={2008},
  publisher={Springer}
}

@article{vo2015data,
  title={Data assimilation and uncertainty assessment for complex geological models using a new PCA-based parameterization},
  author={Vo, Hai X and Durlofsky, Louis J},
  journal={Computational Geosciences},
  volume={19},
  number={4},
  pages={747--767},
  year={2015},
  publisher={Springer}
}

@article{chan2019parametric,
  title={Parametric generation of conditional geological realizations using generative neural networks},
  author={Chan, Shing and Elsheikh, Ahmed H},
  journal={Computational Geosciences},
  volume={23},
  number={5},
  pages={925--952},
  year={2019},
  publisher={Springer}
}

@article{ho2020denoising,
  title={Denoising diffusion probabilistic models},
  author={Ho, Jonathan and Jain, Ajay and Abbeel, Pieter},
  journal={Advances in neural information processing systems},
  volume={33},
  pages={6840--6851},
  year={2020}
}

@article{song2020denoising,
  title={Denoising diffusion implicit models},
  author={Song, Jiaming and Meng, Chenlin and Ermon, Stefano},
  journal={arXiv preprint arXiv:2010.02502},
  year={2020}
}

@article{wang2023prior,
  title={A prior regularized full waveform inversion using generative diffusion models},
  author={Wang, Fu and Huang, Xinquan and Alkhalifah, Tariq A},
  journal={IEEE transactions on geoscience and remote sensing},
  volume={61},
  pages={1--11},
  year={2023},
  publisher={IEEE}
}

@inproceedings{rombach2022high,
  title={High-resolution image synthesis with latent diffusion models},
  author={Rombach, Robin and Blattmann, Andreas and Lorenz, Dominik and Esser, Patrick and Ommer, Bj{\"o}rn},
  booktitle={Proceedings of the IEEE/CVF conference on computer vision and pattern recognition},
  pages={10684--10695},
  year={2022}
}

@article{wu2023physics,
  title={Physics-embedded inverse analysis with algorithmic differentiation for the earth’s subsurface},
  author={Wu, Hao and Greer, Sarah Y and O’Malley, Daniel},
  journal={Scientific Reports},
  volume={13},
  number={1},
  pages={718},
  year={2023},
  publisher={Nature Publishing Group UK London}
}

@book{patankar2018numerical,
  title={Numerical heat transfer and fluid flow},
  author={Patankar, Suhas},
  year={2018},
  publisher={CRC press}
}

@article{eymard2000finite,
  title={Finite volume methods},
  author={Eymard, Robert and Gallou{\"e}t, Thierry and Herbin, Rapha{\`e}le},
  journal={Handbook of numerical analysis},
  volume={7},
  pages={713--1018},
  year={2000},
  publisher={Elsevier}
}

@article{di2024latent,
  title={Latent diffusion models for parameterization of facies-based geomodels and their use in data assimilation},
  author={Di Federico, Guido and Durlofsky, Louis J},
  journal={Computers \& Geosciences},
  volume={194},
  pages={105755},
  year={2025},
  publisher={Elsevier}
}

@article{di20253d,
  title={3D latent diffusion models for parameterizing and history matching multiscenario facies systems},
  author={Di Federico, Guido and Durlofsky, Louis J},
  journal={arXiv preprint arXiv:2508.16621},
  year={2025}
}

@article{raissi2019physics,
  title={Physics-informed neural networks: A deep learning framework for solving forward and inverse problems involving nonlinear partial differential equations},
  author={Raissi, Maziar and Perdikaris, Paris and Karniadakis, George E},
  journal={Journal of Computational physics},
  volume={378},
  pages={686--707},
  year={2019},
  publisher={Elsevier}
}

@book{versteeg2007introduction,
  title={An introduction to computational fluid dynamics the finite volume method, 2/E},
  author={Versteeg, Henk Kaarle},
  year={2007},
  publisher={Pearson Education India}
}

@article{aavatsmark2002introduction,
  title={An introduction to multipoint flux approximations for quadrilateral grids},
  author={Aavatsmark, Ivar},
  journal={Computational Geosciences},
  volume={6},
  number={3},
  pages={405--432},
  year={2002},
  publisher={Springer}
}

@book{hill2007effective,
  title={Effective groundwater model calibration: with analysis of data, sensitivities, predictions, and uncertainty},
  author={Hill, Mary C and Tiedeman, Claire R},
  year={2007},
  publisher={John Wiley \& Sons}
}

@article{lee2014large,
  title={Large-scale hydraulic tomography and joint inversion of head and tracer data using the principal component geostatistical approach (PCGA)},
  author={Lee, Jonghyun and Kitanidis, Peter K},
  journal={Water Resources Research},
  volume={50},
  number={7},
  pages={5410--5427},
  year={2014},
  publisher={Wiley Online Library}
}

@book{Freeze1979,
  title={Groundwater},
  author={Freeze, R. Allan and Cherry, John A.},
  year={1979},
  publisher={Prentice-Hall}
}

@article{Carrera2005,
  title={Computational and conceptual issues in the calibration of seawater intrusion models},
  author={Carrera, Jes{\'u}s and Hidalgo, Juan J and Slooten, Luit J and V{\'a}zquez-Su{\~n}{\'e}, Enric},
  journal={Hydrogeology Journal},
  volume={13},
  pages={204--217},
  year={2005},
  publisher={Springer}
}

@article{Brunner2012,
  title={Hydrogeologic controls on disconnection between surface water and groundwater},
  author={Brunner, Philip and Simmons, Craig T},
  journal={Water Resources Research},
  volume={48},
  number={1},
  year={2012},
  publisher={Wiley Online Library}
}

@article{Yeh2000,
  title={Hydraulic tomography: a new aquifer test method for reconstructing three-dimensional hydraulic conductivity distributions},
  author={Yeh, T-C Jim and Liu, Shlomo},
  journal={Water resources research},
  volume={36},
  number={9},
  pages={2517--2532},
  year={2000},
  publisher={Wiley Online Library}
}

@article{emerick2012ensemble,
  title={Ensemble smoother with multiple data assimilation},
  author={Emerick, Alexandre A and Reynolds, Albert C},
  journal={Computers \& geosciences},
  volume={55},
  pages={3--15},
  year={2012},
  publisher={Elsevier}
}

@book{hirsch2007numerical,
  title={Numerical computation of internal and external flows: The fundamentals of computational fluid dynamics},
  author={Hirsch, Charles},
  year={2007},
  publisher={Elsevier}
}

@article{hinton2006reducing,
  title={Reducing the dimensionality of data with neural networks},
  author={Hinton, Geoffrey E and Salakhutdinov, Ruslan R},
  journal={science},
  volume={313},
  number={5786},
  pages={504--507},
  year={2006},
  publisher={American Association for the Advancement of Science}
}

@inproceedings{masci2011stacked,
  title={Stacked convolutional auto-encoders for hierarchical feature extraction},
  author={Masci, Jonathan and Meier, Ueli and Cire{\c{s}}an, Dan and Schmidhuber, J{\"u}rgen},
  booktitle={International conference on artificial neural networks},
  pages={52--59},
  year={2011},
  organization={Springer}
}

@article{zhan2025toward,
  title={Toward artificial general intelligence in hydrogeological modeling with an integrated latent diffusion framework},
  author={Zhan, Chuanjun and Dai, Zhenxue and Jiao, Jiu Jimmy and Soltanian, Mohamad Reza and Yin, Huichao and Carroll, Kenneth C},
  journal={Geophysical Research Letters},
  volume={52},
  number={3},
  pages={e2024GL114298},
  year={2025},
  publisher={Wiley Online Library}
}

@article{lee2025latent,
  title={Latent diffusion model for conditional reservoir facies generation},
  author={Lee, Daesoo and Ovanger, Oscar and Eidsvik, Jo and Aune, Erlend and Skauvold, Jacob and Hauge, Ragnar},
  journal={Computers \& Geosciences},
  volume={194},
  pages={105750},
  year={2025},
  publisher={Elsevier}
}

@article{shinozuka1996simulation,
  title={Simulation of multi-dimensional Gaussian stochastic fields by spectral representation},
  author={Shinozuka, Masanobu and Deodatis, George},
  journal={Applied Mechanics Reviews},
  volume={49},
  number={1},
  pages={29--53},
  year={1996},
  publisher={American Society of Mechanical Engineers Digital Collection}
}

@article{zhang2024non,
  title={Non-gaussian hydraulic conductivity and potential contaminant source identification: A comparison of two advanced DLPM-based inversion framework},
  author={Zhang, Xun and Jiang, Simin and Wei, Junze and Wu, Chao and Xia, Xuemin and Wang, Xinshu and Zheng, Na and Xing, Jingwen},
  journal={Journal of Hydrology},
  volume={638},
  pages={131540},
  year={2024},
  publisher={Elsevier}
}

@inproceedings{ronneberger2015u,
  title={U-net: Convolutional networks for biomedical image segmentation},
  author={Ronneberger, Olaf and Fischer, Philipp and Brox, Thomas},
  booktitle={International Conference on Medical image computing and computer-assisted intervention},
  pages={234--241},
  year={2015},
  organization={Springer}
}

@inproceedings{frostig2019compiling,
  title={Compiling machine learning programs via high-level tracing},
  author={Frostig, Roy and Johnson, Matthew James and Leary, Chris},
  booktitle={SysML conference 2018},
  year={2019}
}

@article{jax2018github,
  title={JAX: composable transformations of Python+ NumPy programs},
  author={Bradbury, James and Frostig, Roy and Hawkins, Peter and Johnson, Matthew James and Leary, Chris and Maclaurin, Dougal and Necula, George and Paszke, Adam and VanderPlas, Jake and Wanderman-Milne, Skye and others},
  year={2018}
}

@article{deutsch1992geostatistical,
  title={Geostatistical software library and user’s guide},
  author={Deutsch, Clayton V and Journel, Andre G and others},
  journal={New York},
  volume={119},
  number={147},
  pages={578},
  year={1992}
}

@article{cui2014likelihood,
  title={Likelihood-informed dimension reduction for nonlinear inverse problems},
  author={Cui, Tiangang and Martin, James and Marzouk, Youssef M and Solonen, Antti and Spantini, Alessio},
  journal={Inverse Problems},
  volume={30},
  number={11},
  pages={114015},
  year={2014},
  publisher={IOP Publishing}
}

@article{yeung2024gaussian,
  title={Gaussian process regression and conditional Karhunen-Lo{\`e}ve models for data assimilation in inverse problems},
  author={Yeung, Yu-Hong and Barajas-Solano, David A and Tartakovsky, Alexandre M},
  journal={Journal of Computational Physics},
  volume={502},
  pages={112788},
  year={2024},
  publisher={Elsevier}
}

@article{tartakovsky2021physics,
  title={Physics-informed machine learning with conditional Karhunen-Lo{\`e}ve expansions},
  author={Tartakovsky, Alexandre M and Barajas-Solano, David A and He, Qizhi},
  journal={Journal of Computational Physics},
  volume={426},
  pages={109904},
  year={2021},
  publisher={Elsevier}
}

@article{gelbrecht2023differentiable,
  title={Differentiable programming for Earth system modeling},
  author={Gelbrecht, Maximilian and White, Alistair and Bathiany, Sebastian and Boers, Niklas},
  journal={Geoscientific Model Development},
  volume={16},
  number={11},
  pages={3123--3135},
  year={2023},
  publisher={Copernicus Publications G{\"o}ttingen, Germany}
}

@article{shen2023differentiable,
  title={Differentiable modelling to unify machine learning and physical models for geosciences},
  author={Shen, Chaopeng and Appling, Alison P and Gentine, Pierre and Bandai, Toshiyuki and Gupta, Hoshin and Tartakovsky, Alexandre and Baity-Jesi, Marco and Fenicia, Fabrizio and Kifer, Daniel and Li, Li and others},
  journal={Nature Reviews Earth \& Environment},
  volume={4},
  number={8},
  pages={552--567},
  year={2023},
  publisher={Nature Publishing Group UK London}
}

@article{du2025jax,
  title={JAX-MPM: A Learning-Augmented Differentiable Meshfree Framework for GPU-Accelerated Lagrangian Simulation and Geophysical Inverse Modeling},
  author={Du, Honghui and He, QiZhi},
  journal={arXiv preprint arXiv:2507.04192},
  year={2025}
}

@article{du2024neural,
  title={Neural-Integrated Meshfree (NIM) Method: A differentiable programming-based hybrid solver for computational mechanics},
  author={Du, Honghui and He, QiZhi},
  journal={Computer Methods in Applied Mechanics and Engineering},
  volume={427},
  pages={117024},
  year={2024},
  publisher={Elsevier}
}

@article{innes2019differentiable,
  title={A differentiable programming system to bridge machine learning and scientific computing},
  author={Innes, Mike and Edelman, Alan and Fischer, Keno and Rackauckas, Chris and Saba, Elliot and Shah, Viral B and Tebbutt, Will},
  journal={arXiv preprint arXiv:1907.07587},
  year={2019}
}

@article{baydin2018automatic,
  title={Automatic differentiation in machine learning: a survey},
  author={Baydin, Atilim Gunes and Pearlmutter, Barak A and Radul, Alexey Andreyevich and Siskind, Jeffrey Mark},
  journal={Journal of machine learning research},
  volume={18},
  number={153},
  pages={1--43},
  year={2018}
}

@article{chen2023fracture,
  title={Fracture network characterization with deep generative model based stochastic inversion},
  author={Chen, Guodong and Luo, Xin and Jiao, Jiu Jimmy and Jiang, Chuanyin},
  journal={Energy},
  volume={273},
  pages={127302},
  year={2023},
  publisher={Elsevier}
}

@article{sobel19683x3,
  title={A 3x3 isotropic gradient operator for image processing},
  author={Sobel, Irwin and Feldman, Gary and others},
  journal={a talk at the Stanford Artificial Project},
  volume={1968},
  pages={271--272},
  year={1968}
}

@article{golmohammadi2015group,
  title={Group-sparsity regularization for ill-posed subsurface flow inverse problems},
  author={Golmohammadi, Azarang and Khaninezhad, Mohammad-Reza M and Jafarpour, Behnam},
  journal={Water Resources Research},
  volume={51},
  number={10},
  pages={8607--8626},
  year={2015},
  publisher={Wiley Online Library}
}

@article{dwight2006effect,
  title={Effect of approximations of the discrete adjoint on gradient-based optimization},
  author={Dwight, Richard P and Brezillon, Joel},
  journal={AIAA journal},
  volume={44},
  number={12},
  pages={3022--3031},
  year={2006}
}

@article{plessix2006review,
  title={A review of the adjoint-state method for computing the gradient of a functional with geophysical applications},
  author={Plessix, R-E},
  journal={Geophysical Journal International},
  volume={167},
  number={2},
  pages={495--503},
  year={2006},
  publisher={Blackwell Publishing Ltd Oxford, UK}
}
\end{document}